\documentclass[12pt,a4paper]{article}

\usepackage{fullpage}
\usepackage{caption}
\usepackage{subcaption}
\usepackage{amssymb,amsmath,latexsym,amsthm}
\usepackage[english]{babel}
\usepackage{graphicx}
\usepackage{enumerate}
\usepackage{mathdots}
\usepackage{mathrsfs,amsfonts,dsfont}
\usepackage[toc,page]{appendix}
\usepackage{color}
\usepackage[colorlinks=true, linkcolor=black]{hyperref}
\usepackage{tikz-cd}
\usepackage{mathtools}
\usepackage{lscape}
\usepackage{makecell}

\allowdisplaybreaks

\newtheorem{theorem}{Theorem}[section]
\newtheorem{definition}[theorem]{Definition}
\newtheorem{proposition}[theorem]{Proposition}

\newtheorem{examplex}[theorem]{Example}
\newenvironment{example}
{\pushQED{\qed}\examplex}
{\popQED\endexamplex}

\newtheorem{remark}[theorem]{Remark}

\newtheorem{lemma}[theorem]{Lemma}
\newtheorem*{proof*}{Proof}

\newcommand{\extd}{{\rm d}}
\newcommand{\id}{{\rm id}}

\renewenvironment{abstract}
{
\small
\begin{center}
\bfseries \abstractname\vspace{-.5em}\vspace{0pt}
\end{center}
\list{}{%
\setlength{\leftmargin}{5mm}
\setlength{\rightmargin}{\leftmargin}%
}%
\item\relax}
{\endlist
}

\title{Noncommutative Fibre Bundles via Bimodules}
\author{Edwin J. Beggs \and James E. Blake \and \\ Swansea University Department of Mathematics}
\date{}

\begin{document}
\pagenumbering{arabic}

\linespread{1.0}

\maketitle

\begin{abstract}
We construct a Leray-Serre spectral sequence for fibre bundles for de Rham sheaf cohomology on noncommutative algebras.
The morphisms are bimodules with zero-curvature extendable bimodule connections. This generalises definitions involving differentiable algebra maps to differentiable completely positive maps by using the KSGNS construction and Hilbert $C^*$-bimodules with bimodule connections.
We give examples of noncommutative fibre bundles, involving group algebras, matrix algebras, and the quantum torus.
%
%
\end{abstract}

\section{Introduction}

In classical topology, locally trivial fibre bundles capture the idea of a space being locally but not globally a product space.
Well known examples of classical fibre bundles include the M\"{o}bius strip, the Hopf fibration on $SU_2$ and gauge bundles in physics.

A fibre bundle is a map $\pi : E \to B$ from the total space to the base space, satisfying the property that there exists a third space $F$ called the fibre which can be associated in a continuous manner with the pre-image $\pi^{-1}\{b\}$ for each $b \in B$.
For each fibre bundle there is a Leray-Serre spectral sequence, which allows calculation of the cohomology of the total space with coefficients in a sheaf.
For a detailed reference on classical fibre bundles and spectral sequences, see chapters 2 and 9 of \cite{Spanier} respectively.

\indent \quad There are several approaches to generalising the theory of fibre bundles to noncommutative geometry.
In all cases, $B$ is taken to be an algebra of functions on a hypothetical base space, and $A$ as the algebra corresponding to the total space.
Since switching from spaces to algebras reverses arrows, a noncommutative fibre bundle goes from $B$ to $A$.

\indent\quad The noncommutative principal bundle definition using Hopf-Galois extensions has many uses, but is limited to examples with a high degree of symmetry, and cannot provide a general framework corresponding to classical fibre bundles. A definition generalising this is given in \cite{Kaygun19}, this requires a fibre algebra (standing for the Hopf algebra in Hopf-Galois extensions) and an order reversal map (called a left transposition in \cite{Kaygun19}) and  a resulting algebra extension.
In \cite{Echterhoff-fibrations}, there is a definition of noncommutative fibre bundles giving rise to a Leray-Serre spectral sequence, but is limited to the case where the base space algebra is the algebra of functions on a locally compact Hausdorff space.
Another approach is \cite{BegMasLeray}, which builds on \cite{bbsheaf}, defining a noncommutative differential fibre bundle as an algebra map between two algebras satisfying a condition on the quotients of the calculi, and constructs a Leray-Serre spectral sequence converging to the de Rham sheaf cohomology of the total space with coefficients in a module with a zero-curvature connection.
The noncommutative Hopf principal bundle on quantum $SU_2$ is an example of both a Hopf-Galois extension and a differential fibre bundle in this sense.

\indent\quad
For $C^*$ algebras it is known that there are not enough algebra maps
to produce an analogue of classical homotopy theory, and the idea of algebra maps is generalised by asymptotic morphisms \cite{conHig}, $KK$-theory \cite{KaspKK}  and (most relevant from our point of view) to completely positive maps via the KSGNS construction \cite{Lance} using Hilbert $C^*$-modules.
In this paper, we reformulate the definition \cite{BegMasLeray} of differential noncommutative fibre bundles so as to no longer require an algebra map but rather a bimodule with connection.
This is part of the the data required to give a differentiable completely positive map
(Prop.\ 4.86 in \cite{QRG}). We use this to give a Hilbert $C^*$-bimodule example of the noncommutative torus over the circle in
Section~\ref{fibre-bundles-example-circle-in-torus}. Note that specifying a completely positive map requires a choice of element of the Hilbert $C^*$-bimodule, but the spectral sequence does not depend on this choice.


The outline of this paper is as follows.
We begin with a review of spectral sequences and the framework of noncommutative differential geometry that we use throughout the paper, before making our key definition of a noncommutative differential fibre bundle.
The key nontrivial idea that motivates this paper is that when expressing differential forms on a fibre bundle of degree $p$ in the fibre and $q$ in the base, the bimodule map $\sigma_E$ associated with a bimodule connection on a bimodule $E$ can be used instead of an algebra map.
Using this we can define a filtration, and we end up getting a Leray-Serre spectral sequence, except now it works for examples where there isn't an algebra map.
We conclude by calculating two finite-dimensional examples of bimodule differential fibre bundles, followed by one infinite-dimensional example.

Our first example of a bimodule differential fibre bundle is between group algebras, with base algebra $\mathbb{C}G$ and total algebra $\mathbb{C}X$, for a subgroup $G \subset X$.
This happens to come from a differentiable algebra map, and so could also have been calculated using existing theory, but since it can be nicely calculated in full it serves to illustrate our theory.

Our second example is between complex-valued matrix algebras, with base space algebra $M_2(\mathbb{C})$ and total space algebra $M_3(\mathbb{C})$.
The bimodule from this example gives a differentiable map which is completely positive but not an algebra map.

Our third example has base space algebra the circle algebra $\mathbb{C}[S^1]$ and total space algebra the quantum torus $\mathbb{C}_{\theta}[\mathbb{T}^2]$.
This example is infinite dimensional, and also does not come from an algebra map.

This paper is based on one of the chapters from the second author's PhD thesis \cite{Blake24}.

\section{Background}
\subsection{Spectral Sequences}

\begin{figure}[!ht]
\centering
\includegraphics[width=0.8\linewidth]{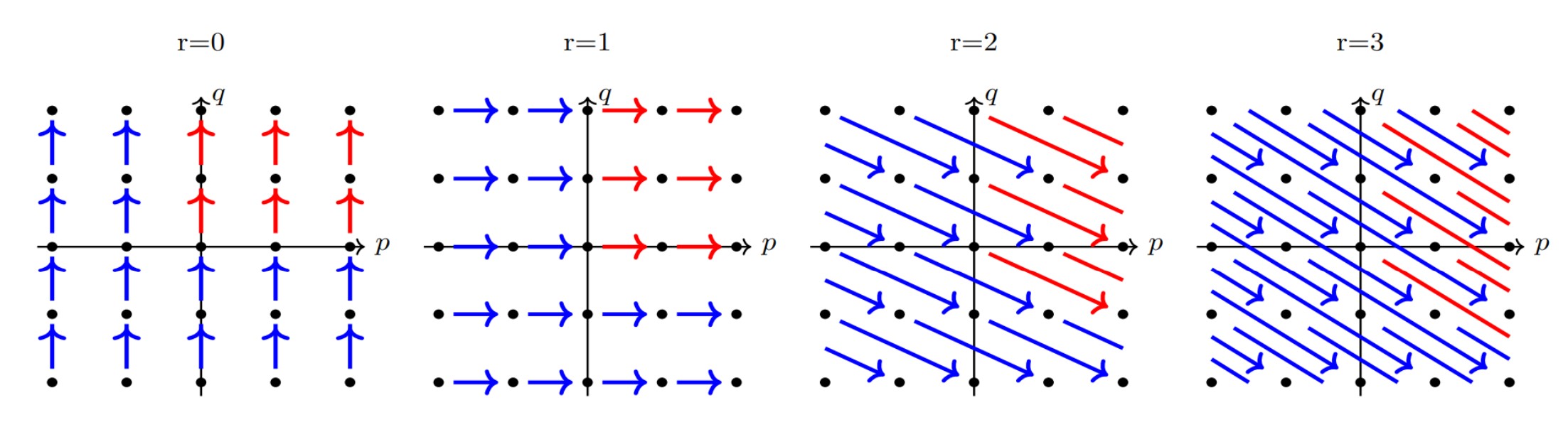}
\caption{Illustration of successive pages of a spectral sequence \cite{Spectral Sequence Image Wikipedia}}
\label{fig:wikipedia---spectralsequencevisualization}
\end{figure}

A spectral sequence $(E_r, d_r)$ is a series of two-dimensional lattices called pages, denoted $E_r$, with each page number $r$ having entry $E^{p,q}_r$ in position $(p,q) \in \mathbb{Z}^2$ and differentials $\extd_r : E^{p,q}_r \to E^{p+r,q+1-r}_r$ satisfying $\extd_r^2 = 0$.
By convention, the $p$-axis is horizontal and the $q$-axis is vertical.
In the case we consider, only the top-right quadrant ($p,q \geq 0$) of page 0 has nonzero entries.
The differentials $\extd_r$ on the $r$th page move right by $r$ entries and down by $r-1$.
Seeing as the differentials square to zero, we can take their cohomology.
The $(r+1)$th page is defined as the cohomology of the $r$th page.
Figure \ref{fig:wikipedia---spectralsequencevisualization} illustrates what a spectral sequence looks like on pages $r = 0,1,2,3$.

A spectral sequence converges if there is a fixed page after which all subsequent pages are the same.
Once taking cohomology no longer changes a spectral sequence, it is said to have stabilised, and position $(p,q)$ on the stable pages is denoted as $E^{p,q}_{\infty}$.

The spectral sequence we use is a variant of the Leray-Serre spectral sequence, which arises from a filtration \cite{McCleary}.
\begin{definition}
Given a cochain complex $C^n$ of vector spaces with linear differential $\extd : C^n \to C^{n+1}$ satisfying $\extd^2 = 0$, we say that a sequence of subspaces $F^mC \subset C$ for $m \geq 0$ is a \textit{decreasing filtration} of $C$ if the following three conditions are satisfied.
\newline\indent 1)\quad 
$\extd{F^mC} \subset F^mC$ for all $m \geq 0$.
\newline\indent 2)\quad 
$F^{m+1}C \subset F^mC$ for all $m \geq 0$.
\newline\indent 3)\quad 
$F^0C = C$ and $F^mC^n := F^mC \cap C^n = \{0\}$ for all $m > n$.
\end{definition}
Given such a filtration, the spectral sequence with first page $E^{p,q}_1 = H^{p+q}(\frac{F^pC}{F^{p+1}C})$ converges to $H^*(C,\extd)$ in the sense that $H^k(C,\extd) = \oplus_{p+q=k} E^{p,q}_{\infty}$.
This can be read off the stabilised sequence as the direct sum along the north-west to south-east diagonals.

\subsection{Bimodules and Connections}

The idea of a bimodule connection was introduced in \cite{DVMic}, \cite{DVMass} and \cite{Mourad} and used in \cite{FioMad}, \cite{Madore}. It was used to construct connections on tensor products in \cite{BDMS}.

\begin{definition}
A first order calculus on an (associative) algebra $A$ is an $A$-bimodule $\Omega^1_A$ satisfying $\Omega^1_A = \text{span}\{a'\extd{a} \mid a,a' \in A \}$ for some linear map $\extd : A \to \Omega^1_A$ with $\extd(ab) = a\,  \extd{b} + \extd{a}\, b$.
\end{definition}
We call the calculus connected if $\ker{\extd} = \mathbb{K}.1$, where $\mathbb{K}$ is the field of scalars of $A$ (in all our examples we take $\mathbb{K} = \mathbb{C}$).
There is a linear map $\wedge : \Omega^p \otimes \Omega^q \to \Omega^{p+q}$, and every element of the higher calculi $\Omega^n$ is a wedge product of elements of $\Omega^1$.
The map $\extd$, called an exterior derivative, can be extended to $\extd: \Omega^n \to \Omega^{n+1}$ by $\extd(\xi \wedge \eta) = \extd{\xi} \wedge \eta + (-1)^{|\xi|}\xi \wedge \extd \eta$ and satisfies $\extd^2 = 0$.

\begin{definition}
Suppose $B$ and $A$ are algebras with calculi $\Omega^1_B$ and $\Omega^1_A$ respectively.
A right bimodule connection $(\nabla_E, \sigma_E)$ on a $B$-$A$ bimodule $E$ consists of two components.

\textbf{(1)} A linear map $\nabla_E : E \to E \otimes_A \Omega^1_A$ called a right connection, satisfying $\nabla_E(e.b) = \nabla_E(e).b + e \otimes \extd{b}$ for $e \in E$, $b \in B$.

\textbf{(2)} A bimodule map $\sigma_E : \Omega^1_B \otimes_B E \to E \otimes_A \Omega^1_A$ called a generalised braiding, satisfying $\nabla_E(be) = \sigma_E(\extd{b} \otimes e) + b \nabla_E(e)$ for $e \in E$, $b \in B$.
\end{definition}

A right connection $\nabla_E$ has curvature $R_E = (\id \otimes \extd + \nabla_E \wedge \id) \nabla_E : E \to E \otimes_A \Omega^{2}_A$, which is always a right module map.
A connection is called flat if its curvature is zero.

The bimodule map $\sigma_E$ is said to be extendable (as in Definition 4.10 of \cite{QRG}) to higher calculi if it extends to a map $\sigma_E : \Omega^n_B \otimes_B E \to E \otimes_A \Omega^{n}_A$ by the formula $\sigma_E(\xi \wedge \eta \otimes e) = (\sigma_E \wedge \id)(\xi \otimes \sigma_E(\eta \otimes e))$.

\begin{lemma}
(See Lemma 1.32 of \cite{QRG})
Every first order calculus $\Omega^1$ on $A$ has a 'maximal prolongation' $\Omega_{max}$ to an exterior algebra, where for every relation $\sum_i a_i.\extd{b_i} = \sum_j \extd{r_j}.s_j$ in $\Omega^1$ for $a_i, b_i, r_j, s_j \in A$ we impose the relation $\sum_i \extd{a_i} \wedge \extd{b_i} = - \sum_j \extd{r_j} \wedge \extd{s_j} \in \Omega^2_{\max}$.
This is extended to higher forms, but no new relations are added.
\end{lemma}

\begin{proposition}
(See Corollary 5.4 of \cite{extendability-for-bimodules})
If the algebra $A$ has maximal prolongation calculi for its higher calculi, and if the curvature $R_E$ is also a left module map, then extendability of $\sigma_E$ is automatic.
\end{proposition}

Next, there are a few results which are proven in \cite{QRG} for left connections but which we need to state for right connections, since the right handed versions are not necessarily just mirror images of the left handed versions, as we see with the signs in Lemma \ref{fibre bundles-nabla-of-sigma}.
We omit the proofs.
\begin{lemma}
Let $E$ be a $B$-$A$ bimodule with extendable right bimodule connection $(\nabla_E,\sigma_E)$.
The connection $\nabla_E$ extends to higher calculi by defining
\begin{equation}
\nabla_E^{[n]} = \id \otimes \extd + \nabla_E \wedge \id : E \otimes_A \Omega^n_A \to E \otimes_A \Omega^{n+1}_A
\end{equation}
\end{lemma}


\begin{lemma}
Let $E$ be a $B$-$A$ bimodule with extendable right bimodule connection $(\nabla_E,\sigma_E)$.
Then $\nabla_E^{[n+1]} \circ \nabla_E^{[n]} = R_E \wedge \id : E \otimes_A \Omega^n_A \to E \otimes_A \Omega^{n+2}_A$, where $R_E = (\id \otimes \extd + \nabla_E \wedge \id) \nabla_E : E \to E \otimes_A \Omega^{2}_A$ is the curvature of $\nabla_E$.
\end{lemma}

In particular, if $R_E = 0$ then the composition $\nabla_E^{[n+1]} \circ \nabla_E^{[n]} = R_E \wedge \id$ vanishes, making the flat connection $\nabla_E$ a cochain differential.


The following lemma almost mirrors the one on page 304 of \cite{QRG}, although we see that by switching sides a power of $-1$ is introduced.

\begin{lemma} \label{fibre bundles-nabla-of-sigma}
Let $E$ be a $B$-$A$ bimodule with extendable right bimodule connection $(\nabla_E,\sigma_E)$ whose curvature $R_E$ is a left module map.
Then for all $n \geq 1$ we have the following equation:
\begin{equation}
\nabla_E^{[n]} \circ \sigma_E = \sigma_E(\extd \otimes \id) + (-1)^n (\sigma_E \wedge \id)(\id \otimes \nabla_E) : \Omega^n_B \otimes_B E \to E \otimes_A \Omega^{n+1}_A.
\end{equation}
\end{lemma}

\section{Theory: Fibre Bundles (Right-handed Version)}

The paper \cite{BegMasLeray} gives a definition of differential fibre bundles, in which given an algebra map $\pi : B \to A$ which extends to a map $\pi^*$ of differential graded algebras, where differential forms of degree $p$ in the base and $q$ in the fibre are given by the quotient $\frac{\pi^* \Omega^p_B \wedge \Omega^q_A}{\pi^* \Omega^{p+1}_B \wedge \Omega^{q-1}_A}$.
The idea behind this is to take the $(p+q)$-forms which are of at least degree $p$ in the base, and then quotient out forms of degree higher than $p$ in the base.

In this paper we take a similar approach, but instead represent these forms by a quotient that doesn't require an algebra map.
However, it does require a bimodule with bimodule connection.  
Since every *-algebra map is also a completely positive map, and Hilbert C*-bimodules with inner products correspond via the KSGNS construction to completely positive maps, this constitutes a generalisation.
Note that the Leray-Serre spectral sequence depends only on the bimodule and bimodule connections, not the
 completely positive maps.

\begin{proposition} \label{submanifolds-filtration-proposition}
Let $E$ be a $B$-$A$ bimodule with extendable zero-curvature right bimodule connection $(\nabla_E,\sigma_E)$.
For $m \leq n$, the cochain complex $C^n = E \otimes_A \Omega^n_A$ with differential $\extd_C := \nabla_E^{[n]} : C^n \to C^{n+1}$ gives the following filtration.
\begin{equation}
F^mC^n = \mathrm{image} \Big( \sigma_E \wedge \id : \Omega^m_B \otimes_B E \otimes_A \Omega^{n-m}_A \to E \otimes_A \Omega^n_A \Big)
\end{equation}
\end{proposition}

\begin{figure}[!ht]
\centering
\includegraphics[width=0.68\linewidth]{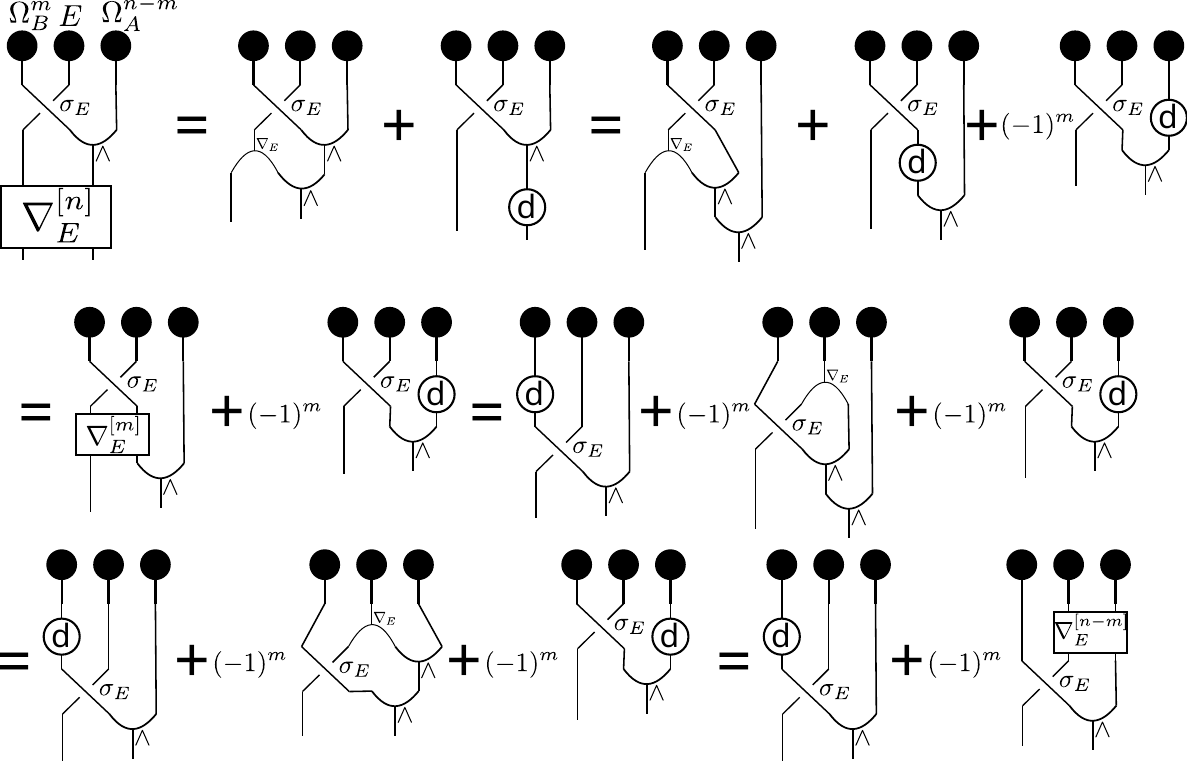}
\caption{Proof that $\extd_C(F^mC) \subset F^mC$}
\label{fig:FB-filtration-differential}
\end{figure}

\begin{figure}[th!]
\centering
\includegraphics[width=0.38\linewidth]{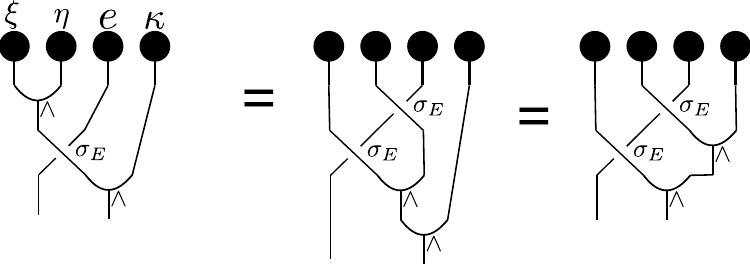}
\caption{Proof the filtration is decreasing}
\label{fig:FB-decreasing-filtration}
\end{figure}

\begin{proof*}
\textbf{(1)} The first property we need for a filtration is $\extd_C{F^mC} \subset F^mC$ for all $m \geq 0$.
This means showing $\nabla_E^{[n]}F^m(E \otimes_A \Omega^n_A) \subset \bigoplus_{n' \geq 0} F^m(E \otimes_A \Omega^{n'}_A)$.

In the calculations in Figure \ref{fig:FB-filtration-differential} we start with $\nabla^{[n]}_E(\sigma_E \wedge \id)$, then use the fact that $\nabla_E^{[n]} = \id \otimes \extd + \nabla_E \wedge \id$, then use associativity of $\wedge$ and expand $\extd \wedge$, then recognise the formula for $\nabla^{[m]}_E$, then use the formula $\nabla_E^{[m]} \circ \sigma_E = \sigma_E(\extd \otimes \id) + (-1)^m (\sigma_E \wedge \id)(\id \otimes \nabla_E)$ we showed earlier, then use associativity of $\wedge$, then recognise the formula for $\nabla^{[n-m]}_E$.

This is in $F^{m+1}C^{n+1} + F^mC^{n+1}$.
However, as we will show in the next step, the filtration is decreasing, so as required, it is contained in $F^mC$.

\textbf{(2)} The second property we need for a filtration is $F^{m+1}C \subset F^mC$ for all $m \geq 0$.

In a differential calculus (as opposed to a more general differential graded algebra), elements of the higher calculi can all be decomposed into wedge products of elements of $\Omega^1$, and so $\Omega^{m+1}_B = \Omega^m_B \wedge \Omega^1_B$.
Let $\xi \in \Omega^m_B$, $\eta \in \Omega^1_B$, $e \in E$, $\kappa \in \Omega^{n-m-1}_A$.
Then $\xi \wedge \eta \otimes e \otimes \kappa \in \Omega^{m+1}_B \otimes_B E \otimes \Omega^{n-m-1}_A$, so the map $\sigma_E \wedge id$ takes it to $E \otimes_A \Omega^n_A$, and the image of all such things is $F^{m+1}C^n$.
We have the string diagram Figure \ref{fig:FB-decreasing-filtration} for $(\id \otimes \wedge)(\sigma \otimes \id)(\wedge \otimes \id \otimes \id)(\xi \otimes \eta \otimes e \otimes \kappa)$, where we use that $\sigma_E$ is extendable and that $\wedge$ is associative.
This shows that $F^{m+1}C^n$ lies in $\textrm{im}(\sigma_E \wedge \id) : \Omega^m_B \otimes_B E \otimes_A \Omega^{n-m}_A \to E \otimes_A \Omega^n_A$, i.e. in $F^mC^n$, and hence that the filtration is decreasing in $m$.

\textbf{(3)} The third property we need is $F^0C = C$.
\begin{align*}
F^0C^n = \textrm{im}(\sigma_E \wedge \id) : B \otimes_B E \otimes_A \Omega^{n}_A \to E \otimes_A \Omega^n_A
\end{align*}
Recalling that $\sigma_E(1 \otimes e) = e \otimes 1$ when $m=0$, the set $F^0C^n$ consists of elements $b.e \otimes \xi$, which gives all of $C^n$.

\textbf{(4)} The final property we need is $F^mC^n := F^mC \cap C^n = \{0\}$ for all $m > n$.
This holds because for $m > n$, we have $\Omega^{n-m} = 0$, giving $F^mC^n = {im(\sigma_E \wedge \id)} : 0 \to C^n$, which has zero intersection with $C^n$.
\qed
\end{proof*}

\begin{definition}
Given a filtration as above, we define \textit{differential forms with coefficients} in $E$ of degree $p$ in the fibre and $q$ in the base as the quotient
\begin{equation} \label{fibre-bundles-M}
M_{p,q} := \frac{F^pC^{p+q}}{F^{p+1}C^{p+q}} = \frac{\sigma_E(\Omega^p_B \otimes_B E) \wedge \Omega^q_A}{\sigma_E(\Omega^{p+1}_B \otimes_B E) \wedge \Omega^{q-1}_A} \ ,
\end{equation}
and from these we denote forms with coefficients in $E$ of degree $q$ in the fibre only as:
\begin{equation} \label{fibre-bundles-N}
N_q := M_{0,q} = \frac{C^q}{F^1C^q} = \frac{E \otimes_A \Omega^q_A}{\sigma_E(\Omega^1_B \otimes_B E) \wedge \Omega^{q-1}_A} \ .
\end{equation}
\end{definition}

\begin{proposition}
Let $E$ be a $B$-$A$ bimodule with extendable zero-curvature right bimodule connection $(\nabla_E,\sigma_E)$.
Then there is a well-defined surjective linear map:
\begin{align} \label{yuop}
g : \Omega^p_B \otimes_B N_q \to M_{p,q} \ ,
\quad g(\xi \otimes [e \otimes \eta]) = [(\sigma_E \wedge \id)(\xi \otimes e \otimes \eta)] \ .
\end{align}
\end{proposition}

\begin{figure}[th!]
\centering
\includegraphics[width=0.3\linewidth]{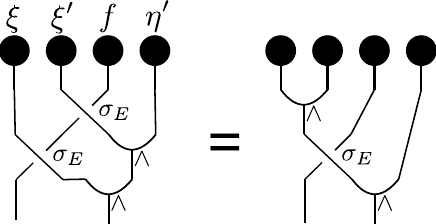}
\caption{Re-writing $g(\xi \otimes [e \otimes \eta])$}
\label{fig:FB-g-isomorphism}
\end{figure}

\begin{proof*}
Surjectivity follows from the definition of the map, so we only need to show that $g$ is well-defined on equivalence classes, i.e. that if $[e \otimes \eta]=0$ then we also have $[(\sigma_E \wedge \id)(\xi \otimes e \otimes \eta)] = 0$.
By definition, we have $[e \otimes \eta] = 0 \in N_q$ if and only if $e \otimes \eta = (\sigma_E \wedge \id)(\xi' \otimes f \otimes \eta')$ for some $\xi' \in \Omega^1_B$, $f \in E$, $\eta' \in \Omega^{q-1}_A$ (summation implicit).
Thus, using associativity of $\wedge$ and then extendability of $\sigma$, we can re-write $g(\xi \otimes [e \otimes \eta])$ as in Figure \ref{fig:FB-g-isomorphism}, which we can see is in the image of $\sigma_E \wedge \id : \Omega^{p+1}_B \otimes_B E \otimes_A \Omega^{q-1}_A \to E \otimes_A \Omega^{p+q}_A$, and hence has equivalence class zero in $M_{p,q}$.
\qed
\end{proof*}

In a classical fibre bundle, the differential forms on the total space would split into a direct sum $\extd x_1 \wedge \cdots \wedge \extd x_p \otimes \extd y_1 \wedge \cdots \wedge \extd y_q$ of forms in the direction of the base space and forms in the direction of the fibre, but in a noncommutative context there is no obvious algebra that can be called the fibre.
Consequently, in the following definition of a bimodule noncommutative fibre bundle (note that the algebra maps approach employs a similar idea), we take the quotient of forms on the total space by forms on the base space as a stand-in for forms on the fibre.
In the classical case where there is a direct sum, this quotient reduces to the usual differential forms on the fibre.

\begin{definition}
For algebras $A$ and $B$ with differential calculi, we call a $B$-$A$ bimodule $E$ a \textit{(bimodule) differential fibre bundle} if it satisfies the following three properties:

\textbf{(1)} There is an extendable zero-curvature right bimodule connection $(\nabla_E,\sigma_E)$ on $E$.

\textbf{(2)} For all $p \geq 1$ the calculi $\Omega^p_B$ are flat as right modules

\textbf{(3)} For all $p,q \geq 0$ the map $g$ in (\ref{yuop}) is an isomorphism.
\end{definition}

\begin{remark}
Recall that flatness of $\Omega^p_B$ as a right $B$-module means that if
\[
\begin{tikzcd}
0 \arrow{r}
&E_1 \arrow{r}{\phi_1}
&E_1 \arrow{r}{\phi_2}
&E_1 \arrow{r}
&0
\end{tikzcd}
\]
is a short exact sequence of left $B$-modules and left $B$-module maps then the following sequence of left $B$-modules and left $B$-module maps is also short exact.
\[
\begin{tikzcd}
0 \arrow{r}
&\Omega^p_B \otimes_B E_1 \arrow{r}{\id \otimes \phi_1}
&\Omega^p_B \otimes_B E_2 \arrow{r}{\id \otimes \phi_2}
&\Omega^p_B \otimes_B E_3 \arrow{r}
&0.
\end{tikzcd}
\]
If $\Omega^p_B$ is finitely generated projective, then the property of being flat is automatic.
\end{remark}

In the remainder of this section we prove the following result corresponding to the classical Leray-Serre spectral sequence.

\begin{theorem}
For algebras $A$ and $B$ with differential calculi, suppose that $(E,\nabla_E,\sigma_E)$ 
is a bimodule differential fibre bundle.
Then there is a spectral sequence for the filtration, which has first page
\begin{align*}
E^{p,q}_1 = H^{p+q}(M_{p,q}) \cong H^{p+q}(\Omega^p_B \otimes_B N_q) = \Omega^p_B \otimes_B H^q(N),
\end{align*}
and second page position $(p,q)$ given by $H^p(B, H^q(N), \nabla_{q})$, and which converges to $H(A,E,\nabla_E)$ in the sense that the entries stabilise at large page number.
\end{theorem}
\begin{proof} 
In Lemma~\ref{lem111} we show that the cohomology giving the first page of the spectral sequence is isomorphic to a construction involving the `fibre' cohomology $H^q(N)$. In Proposition~\ref{prus} we show that $H^q(N)$ is itself a sheaf over the base (i.e.\ that it has a zero curvature connection $\nabla_q$). Thus we get a cochain complex
\[
\begin{tikzcd}
0 \arrow{r}
& H^q(N) \arrow{r}{\nabla_q}
& \Omega^1_B \otimes_B H^q(N) \arrow{r}{\nabla_q^{[1]}}
& \Omega^2_B \otimes_B H^q(N) \arrow{r}{\nabla_q^{[2]}}
& \cdots
\end{tikzcd}
\]
giving the sheaf cohomology $H^p(B,H^q(N),\nabla_q)$, and this is the second page of the 
 spectral sequence.
\end{proof}

The following Lemma corresponds to Lemma 4.64 of \cite{QRG}, but with a number of differences to adapt it from algebra maps fibre bundles to bimodule fibre bundles.

\begin{lemma} \label{lem111}
For a bimodule differential fibre bundle $(E, \nabla_E)$, there is a cochain complex 
\[
\begin{tikzcd}
\cdots \arrow{r}{d}
& M_{p,q-1} \arrow{r}{d}
& M_{p,q} \arrow{r}{d}
& M_{p,q+1} \arrow{r}{d}
& \cdots
\end{tikzcd},
\]
whose differential on $M_{p,q}$ is $[\nabla^{[p+q]}_E]$, and whose cohomology we denote as $H^{q}(M_{p,q})$.
Then by equation (\ref{fibre-bundles-N}), $H^q(N) := H^q(M_{0,q})$.
The isomorphism $g$ is a (graded) cochain map, and extends to the following isomorphism of cohomology:
\begin{equation}
\hat{g} : \Omega^p_B \otimes_B H^q(M_{0,q}) \to H^{q}(M_{p,q}), \qquad
{\xi \otimes [[e \otimes \eta]]} \mapsto {[[\sigma_E(\xi \otimes e) \wedge \eta]]}.
\end{equation}

\end{lemma}

\begin{proof*}
\textbf{(1)} The differential is well-defined by Proposition \ref{submanifolds-filtration-proposition}, and satisfies $\extd^2 = 0$ by flatness of $\nabla_E$.

\textbf{(2)} Recall that $(E, \nabla_E)$ being a bimodule differential fibre bundle implies that for all $p,q \geq 0$ there is an isomorphism $g : \Omega^p_B \otimes_B N_q \to M_{p,q}$.
We need to show that the differential $[\nabla^{[p+q]}_E]$ commutes with $g$, i.e. that the following diagram commutes:
\[
\begin{tikzcd}
\Omega^p_B \otimes_B N_q \arrow{rr}{(-1)^p \id \otimes [\nabla^{[q]}_E]} \arrow{d}{g}
&& \Omega^p_B \otimes_B N_{q+1} \arrow{d}{g} \\
M_{p,q} \arrow{rr}{[\nabla^{[p+q]}_E]}
&& M_{p,q+1}
\end{tikzcd}
\]
In the proof that $(F^mC,\nabla^{[n]}_E)$ is a filtration, we calculated (diagrammatically) that 
\begin{align*}
\nabla^{[p+q]}_E(\sigma_E \wedge \id) = \sigma_E(\extd \otimes \id) \wedge \id + (-1)^p(\sigma_E \wedge \id)(\id \otimes \nabla^{[q]}_E).
\end{align*}
However, we know that the term $\sigma_E(\extd \otimes \id) \wedge \id$ has equivalence class zero in $M_{p,q+1}$.
Taking equivalence classes therefore gives 
\begin{align*}
[\nabla^{[p+q]}_E \circ g] = (-1)^p[(\sigma_E \wedge \id)(\id \otimes \nabla^{[q]}_E)].
\end{align*}
Going the other way around the diagram, we get
\begin{align*}
g((-1)^p \id \otimes [\nabla^{[q]}_E]) = (-1)^p[(\sigma_E \wedge \id)(\id \otimes \nabla^{[q]}_E)].
\end{align*}
These coincide, so the diagram commutes.

\textbf{(3)} Secondly, we need to show that $g$ extends to cohomology, i.e. that the map
\begin{align*}
\Omega^p_B \otimes_B H^q(N_{q}) \to H^{q}(M_{p,q}), \qquad
{\xi \otimes [[e \otimes \eta]]} \mapsto {[[\sigma_E(\xi \otimes e) \wedge \eta]]}
\end{align*}
is an isomorphism.
Make the following two definitions.
\begin{align*}
Z_{p,q} := \rm{im}(\extd) : M_{p,q-1} \to M_{p,q}, \qquad
K_{p,q} := \ker(\extd) : M_{p,q} \to M_{p,q+1}
\end{align*}
Therefore $H^{p+q}(M_{p,q}) = \frac{K_{p,q}}{Z_{p,q}}$.

Next we show that the differential $\extd : M_{0,q} \to M_{0,q+1}$ is a left $B$-module map.
We take $[\nabla^{[q]}_E(b.e \otimes \eta)]$ and apply the definition of $\nabla^{[q]}$, then use the Leibniz rule, then use the fact that $[\sigma_E(\extd{b} \otimes e) \wedge \eta] = 0$ to calculate:
\begin{align*}
[\nabla^{[q]}_E & (b.e \otimes \eta)]
= [b.e \otimes \extd{\eta}] + [\nabla_E(b.e) \wedge \eta] \\
&= [b.e \otimes \extd{\eta}] + [\sigma_E(\extd{b} \otimes e) \wedge \eta] + [b \nabla_E(e) \wedge \eta] \\
&= [b(\id \otimes \extd + \nabla_E \wedge \id)(e \otimes \eta)] = [b \nabla^{[q]}_E(e \otimes \eta)].
\end{align*}
Hence there is an exact sequence of left $B$-modules and left $B$-module maps:
\[
\begin{tikzcd}
0 \arrow{r} &K_{0,q} \arrow{r}{\rm{inc}} &M_{0,q} \arrow{r}{\extd} &Z_{0,q+1} \arrow{r} &0
\end{tikzcd}
\]
Taking the tensor product with the flat right module $\Omega^p_B$ gives another exact sequence:
\[
\begin{tikzcd}
0 \arrow{r} &\Omega^p_B \otimes_B K_{0,q} \arrow{r}{\id \otimes \rm{inc}} &\Omega^p_B \otimes_B M_{0,q} \arrow{r}{\id \otimes \extd} &\Omega^p_B \otimes_B Z_{0,q+1} \arrow{r} &0
\end{tikzcd}
\]
Applying $g$ to the elements of this sequence, the first part of the proof tells us that the following diagram commutes.
\[
\begin{tikzcd}
0 \arrow{r} &\Omega^p_B \otimes_B K_{0,q} \arrow{r}{\id \otimes \rm{inc}} \arrow{d}{g} &\Omega^p_B \otimes_B M_{0,q} \arrow{r}{\id \otimes \extd} \arrow{d}{g} &\Omega^p_B \otimes_B Z_{0,q+1} \arrow{r} \arrow{d}{g} &0 \\
&M_{p,q} \arrow{r}{\id} & M_{p,q} \arrow{r}{(-1)^p \extd} & M_{p,q+1}
\end{tikzcd}
\]
Note that the middle instance of $g$ is an isomorphism, while the first and third are merely injective.
This diagram gives the following two isomorphisms.
\begin{align*}
Z_{p,q+1} \cong \Omega^p_B \otimes_B Z_{0,q+1}, \qquad
K_{p,q+1} \cong \Omega^p_B \otimes_B K_{0,q+1},
\end{align*}
and consequently
\begin{align*}
Z_{p,q} \cong \Omega^p_B \otimes_B Z_{0,q}, \qquad
K_{p,q} \cong \Omega^p_B \otimes_B K_{0,q}.
\end{align*}
By definition of $H^q(N) = \frac{K_{0,q}}{Z_{0,q}}$, we have another short exact sequence:
\[
\begin{tikzcd}
0 \arrow{r} &Z_{0,q} \arrow{r}{\rm{inc}} &K_{0,q} \arrow{r} &H^q(N) \arrow{r} &0
\end{tikzcd}
\]
Taking the tensor product with the flat right module $\Omega^p_B$ gives the exact sequence:
\[
\begin{tikzcd}
0 \arrow{r} &\Omega^p_B \otimes_B Z_{0,q} \arrow{r}{\id \otimes \rm{inc}} &\Omega^p_B \otimes_B K_{0,q} \arrow{r} &\Omega^p_B \otimes_B H^q(N) \arrow{r} &0
\end{tikzcd}
\]
Therefore 
\begin{align*}
\Omega^p_B \otimes_B H^q(N) \cong \frac{\Omega^p_B \otimes_B K_{0,q}}{\Omega^p_B \otimes_B Z_{0,q}} \cong \frac{K_{p,q}}{Z_{p,q}} = H^{p+q}(M_{p,q}).
\end{align*}
This is the isomorphism we wanted to show.
\qed
\end{proof*}

\begin{proposition} \label{prus}
We have a zero-curvature left connection on the cohomology of the fibre
\begin{equation}
\nabla_{q} : H^q(N) \to \Omega^1_B \otimes_B H^q(N), \qquad
\nabla_{q}([[e \otimes \xi]]) = \eta \otimes [[f \otimes \kappa]]
\end{equation}
where  $\eta \in \Omega^1_B$, $f \in E$, $\kappa \in \Omega^q_A$ are defined by  (summation implicit) 
\begin{align} \label{udre}
\nabla^{[n]}_E(e \otimes \xi) = \sigma_E(\eta \otimes f) \wedge \kappa \in \sigma_E(\Omega^1_B \otimes E) \wedge \Omega^q_A\ .
\end{align}
\end{proposition}

\begin{proof*}
\textbf{(1)} Firstly, we show that the map $\nabla_q$ is well-defined.
Since $H^q(N) = \frac{\ker [\nabla^{[q]}_E]}{\rm{im} [\nabla^{[q-1]}_E]}$ it follows that for all $e \in E$, $\xi \in \Omega^q_A$ such that $[[e \otimes \xi]] \in H^q(N)$, the equivalence class $[\nabla^{[n]}_E(e \otimes \xi)]$ vanishes in $N_{q+1}$.
But for $\nabla^{[n]}_E(e \otimes \xi)$ to lie in the denominator of $N_{q+1} = \frac{E \otimes_A \Omega^{q+1}_A}{\sigma_E(\Omega^1_B \otimes_B E) \wedge \Omega^{q}_A}$ means that there exist $\eta \in \Omega^1_B$, $f \in E$, $\kappa \in \Omega^q_A$satisfying (\ref{udre}).
Applying the isomorphism $\hat{g}^{-1}$ to $[[\nabla^{[n]}_E(e \otimes \xi)]]$ then gives $\eta \otimes [[f \otimes \kappa]]$.

\textbf{(2)} Next, we show for all $b \in B$ that $\nabla_q$ satisfies the left Leibniz rule.
We take $\nabla^{[q]}_E(be \otimes \xi)$ and use the definition of $\nabla^{[q]}_E$ then the braided Leibniz rule for $\nabla_E$ and then once again the definition of $\nabla^{[q]}_E$ to calculate:
\begin{align*}
\nabla^{[q]}_E & (be \otimes \xi) = (\id \otimes \extd + \nabla_E \wedge \id)(be \otimes \xi) 
= be \otimes \extd{\xi} + \nabla_E(be) \wedge \xi \\
&= be \otimes \extd{\xi} + \sigma_E(\extd{b} \otimes e) \wedge \xi + b \nabla_E(e) \wedge \xi \\
&= \sigma_E(\extd{b} \otimes e) \wedge \xi + b \nabla^{[q]}_E(e \otimes \xi)
\end{align*}
Taking equivalence classes and using the isomorphism $\hat{g}$ gives the desired result that
\begin{align*}
\nabla_{q}(b.[[e \otimes \xi]]) = {\extd{b} \otimes [[e \otimes \xi]]} + {b \nabla_{q}([[e \otimes \xi]])}.
\end{align*}

\textbf{(3)} Lastly, we show that the curvature, $R_{q} = {(\extd \otimes \id - \id \wedge \nabla_{q})\nabla_{q}}$ vanishes.

Denoting $\nabla^{[q]}_E(f \otimes \kappa) = \sigma_E(\eta' \otimes f') \wedge \kappa'$, we have:
\begin{align*}
R_q([[e \otimes \xi]]) &= \extd{\eta} \otimes [[f \otimes \kappa]] + \eta \wedge \nabla_q([[f \otimes \kappa]]) \\
&= \extd{\eta} \otimes [[f \otimes \kappa]] + \eta \wedge \eta' \otimes [[f' \otimes \kappa']]
\end{align*}
To show this vanishes, we want to show $\extd{\eta} \otimes [f \otimes \kappa] + \eta \wedge \eta' \otimes [f' \otimes \kappa'] = 0$.

As the curvature $R_E$ of $\nabla_E$ vanishes, we have:
\begin{align*}
0 &= \nabla^{[q+1]}_E \circ \nabla^{[q]}_E(e \otimes \xi) = \nabla^{[q+1]}_E(\sigma_E(\eta \otimes f) \wedge \kappa) \\
&= (\extd \otimes \id + \id \wedge \nabla_E)(\sigma_E(\eta \otimes f) \wedge \kappa)
\end{align*}
Taking equivalence classes and using the isomorphism $g$, we get
\begin{align*}
0 &= (\extd \otimes \id + \id \wedge \nabla_E)(\eta \otimes [f \otimes \kappa]) 
= \extd{\eta} \otimes [f \otimes \kappa] + \eta \wedge \nabla_E([f \otimes \kappa]) \\
&= \extd{\eta} \otimes [f \otimes \kappa] + \eta \wedge \eta' \otimes [f' \otimes \kappa']
\end{align*}
as required.
Hence $R_q = 0$.
\qed
\end{proof*}

\section{Theory: Fibre Bundles (Left-handed Version)}

Swapping left and right, the previous construction can be mirrored to use an $A$-$B$ bimodule $E$ with an extendable zero-curvature left bimodule connection $(\nabla_E, \sigma_E)$, where $\nabla_E : E \to \Omega^1_A \otimes_A E$ and $\sigma_E : E \otimes_B \Omega^1_B \to \Omega^1_A \otimes_A E$.
In this case, zero curvature means $R_E = (\extd \otimes \id - \id \wedge \nabla_E) \nabla_E = 0$.
The bimodule connection satisfies
\begin{align*}
\nabla^{[n]}_E \sigma_E = (\id \wedge \sigma_E)(\nabla_E \otimes \id) + \sigma_E(\id \otimes \id) : E \otimes_B \Omega^n_B \to \Omega^{n+1}_A \otimes_A E.
\end{align*}
The cochain complex $C^n = \Omega^n_A \otimes_A E$ with differential $C^n \to C^{n+1}$ given by
\begin{align*}
\extd_C = \nabla^{[n]}_E = \id \otimes \extd + (-1)^n \nabla_E \wedge \id
\end{align*}
has a filtration
\begin{align*}
F^mC^n = \rm{im}(\id \wedge \sigma_E) : \Omega^{n-m}_A \otimes_A E \otimes_B \Omega^m_B \to \Omega^n_A \otimes_A E.
\end{align*}
The quotients for the fibre are given as follows.
\begin{align*}
&M_{p,q} := \frac{F^pC^{p+q}}{F^{p+1}C^{p+q}} = \frac{\Omega^q_A \wedge \sigma_E(E \otimes_B, \Omega^p_B)}{\Omega^{q-1}_A \wedge \sigma_E(E \otimes_B \Omega^{p+1}_B)} \\
&N_q := M_{0,q} = \frac{C^q}{F^1C^q} = \frac{\Omega^q_A \otimes_A E}{\Omega^{q-1}_A \wedge \sigma_E(E \otimes_B \Omega^1_B)},
\end{align*}
There is then a well-defined map
\begin{align*}
g : N_q \otimes_B \Omega^p_B \to M_{p,q}, \qquad
[\eta \otimes e] \otimes \xi \mapsto [(\id \wedge \sigma)(\eta \otimes e \otimes \xi)]
\end{align*}
which extends to cohomology.

We say that $E$ is a differential fibre bundle if $g$ is an isomorphism for all $p,q \geq 0$ and if the calculi $\Omega^p_B$ are flat as left modules for all $p \geq 0$

On the cohomology we have the following a zero-curvature right connection.
\begin{align*}
\nabla_q : H^q(N) \to H^q(N) \otimes_B \Omega^1_B,
&& \nabla_q([[\xi \otimes e]]) = [[\kappa \otimes f]] \otimes \eta,
\end{align*}
where $\nabla^{[q]}_E(\xi \otimes e) = \kappa \wedge \sigma_E(f \otimes \eta) \in \Omega^q_A \wedge \sigma_E(E \otimes_B \Omega^1_B) \subset \Omega^{q+1}_A \otimes_A E$ with summation implicit.
Assuming that we have a differential fibre bundle, there is then a spectral sequence converging to $H(A,E,\nabla_E)$ with first page position $(p,q)$ given by $E^{p,q}_1 = H^q(N) \otimes_B \Omega^p_B$ and second page position $(p,q)$ given by $H^p(B, H^q(N), \nabla_q)$.

\section{Example: Group Algebras} \label{fibre-bundles-example-group-algebras}

\begin{example} \label{algebra-index-group-algebra} \cite{Quantum-groups-calculus}
For a finite group $X$, the group algebra $\mathbb{C}X$ is the linear extension of the group.
It has star structure $x^* = x^{-1}$ for all $x \in X$.
Its left translation-invariant calculi are given by right $\mathbb{C}X$-modules $V$ and maps $\omega : X \to V$ satisfying
the cocycle condition
\begin{align*}
\omega(xy) = \omega(x) \triangleleft y + \omega(y), \quad \forall x,y \in X\ .
\end{align*}
The calculus is then given as a free module $\Omega^1 = \mathbb{C}X.V$.
Each pair $(V, \omega)$ of a right $\mathbb{C}X$-module and a cocycle therefore gives a translation-invariant calculus.
The right action on the calculus is $v.x = x.(v \triangleleft x)$, and the differential is $\extd{x} = x \, \omega(x)$, for all $x \in X$ and $v \in V$.
The calculus is connected if and only if $\omega(x) \neq 0$ for all $x \in X \backslash \{ e \}$.
\end{example}

Note:  It follows that $\omega(x^{-1}) = - \omega(x) \triangleleft x^{-1}$ and that $\omega(1) = 0$ ($1$ is the identity of the group and group algebra). 
We can also write $\Omega^1 = V . \mathbb{C}X$ by using $v.x = x. (v \triangleleft x)$. A usual notation for left translation-invariant first order calculi on Hopf algebras is $\Lambda^1$, so we could refer to $V$ above as $\Lambda^1_{\mathbb{C} X}$.

\begin{lemma}
Let $X$ be a finite group with calculus given by a right representation $V$ and a cocycle $\omega : \mathbb{C}X \to V$, and which has a subgroup $G$.
Then subspace $W$ of $V$ spanned by $\omega(g)$ for all $g \in G$ is a right representation of $G$, and has complement $W^{\perp}$ which is also right representation of $G$. The restriction of $\omega$ to a cocycle $\mathbb{C}G \to W$ gives a calculus on the subgroup $G$.
\end{lemma}
\begin{proof*}
The cocycle condition $\omega(x) \triangleleft y = \omega(xy) - \omega(y)$ defines a right action on $W$, which gives a calculus on $\mathbb{C}G$.
Since $G$ is a finite group, the representation $V$ has an invariant inner product $\overline{V} \otimes V \to C$ (invariant meaning $\langle \overline{v \triangleleft g} , v \triangleleft g \rangle = \langle \overline{v} , v \rangle$), from which it follows that $V = W \oplus W^{\perp}$, where $W^{\perp}$ is the perpendicular complement of $W$.
The vector space $W^{\perp}$ is then also a representation of $G$.
\qed
\end{proof*}

Write $\pi^{\perp}:V \to W^{\perp}$ and $\pi:V \to W$ for the projections corresponding to the direct sum.

\begin{proposition}
If for the higher calculi on $\mathbb{C}X$ we assume that $\extd V= 0$, then the wedge product $\wedge$ is antisymmetric on
elements of $V$, thus the higher left invariant forms $\Lambda_{\mathbb{C} X}$ can be taken to be the antisymmetric tensor algebra on $V$ (the maximal prolongation calculus).
\end{proposition}

\begin{proof*}
Since $v \triangleleft x = x^{-1} v x$, it follows that $x(v \triangleleft x) = vx$.
Applying $\extd$ to this and using the assumption that $\extd(V) = 0$, we obtain $\extd{x} \wedge (v \triangleleft x) = -v \wedge \extd{x}$.
Using this,
\begin{align*}
v \wedge \omega(x) &= v \wedge (x^{-1} \extd{x})
= (vx^{-1}) \wedge \extd{x}
= x^{-1}(v \triangleleft x^{-1}) \wedge \extd{x}
= -x^{-1} \extd{x} \wedge v
= - \omega(x) \wedge v.
\end{align*}
Since the images $\omega(x)$ span $V$, this proves that $\wedge$ is antisymmetric on $\Lambda_{\mathbb{C}X}$.
\qed
\end{proof*}

Now we look at fibre bundles.
Suppose $G$ is a finite subgroup of a group $X$, and take $A = \mathbb{C}X$, $B = \mathbb{C}G$ as in the discussion of fibre bundles earlier.
Equip $\mathbb{C}X$ with calculus as above for $\Lambda^1_{\mathbb{C}X} = V$ and some cocycle $\omega : \mathbb{C}X \to V$ for some right representation $V$ of $\mathbb{C}X$.
For the higher calculi on $\mathbb{C}X$ take maximal prolongation plus the assumption $\extd V=0$.
For the calculus on $\mathbb{C}G$ take $\Lambda^1_{\mathbb{C}G} = W = \omega(\mathbb{C}G)$ with cocycle the restriction of $\omega$ to $\mathbb{C}G$, and maximal prolongation for the higher calculi.

\begin{proposition} Take $E = \mathbb{C}X$ as a
 $\mathbb{C}G$-$\mathbb{C}X$ bimodule  with left and right actions given by multiplication.
When the algebras are equipped with the calculi above there is a zero-curvature extendable right bimodule connection on $E$ given by $(\nabla_E, \sigma_E)$, where
\begin{align*}
\nabla_E &: \mathbb{C}X \to \mathbb{C}X \otimes_{\mathbb{C}X} \Omega^1_{\mathbb{C}X}, \qquad
x \mapsto 1 \otimes \extd{x}\ ,\cr
\sigma_E &: \Omega^1_{\mathbb{C}G} \otimes_{\mathbb{C}G} \mathbb{C}X \to \mathbb{C}X \otimes_{\mathbb{C}X} \Omega^1_{\mathbb{C}X}, \qquad
\extd{g} \otimes x \mapsto 1 \otimes \extd{g}.x.
\end{align*}
\end{proposition}
\begin{proof*}
The connection satisfies the condition $\nabla_E(gx) = \sigma_E(\extd{g} \otimes x) + g \nabla_E(x)$ required to be a bimodule connection, since $\sigma_E(\extd{g} \otimes x) = 1 \otimes (\extd(gx) - g \extd{x}) = 1 \otimes \extd{g}.x$.
The curvature is zero because $\extd$ has zero curvature.
The connection is extendable as $\sigma_E(\xi \otimes x) = 1 \otimes \xi.x$ for all $\xi \in \Omega^n_{\mathbb{C}G}$.
\qed
\end{proof*}

\begin{proposition}
Equip $A = \mathbb{C}X$, $B = \mathbb{C}G$ with calculi as above, the $B$-$A$ bimodule $E = \mathbb{C}X$ with actions given by multiplication.
The right bimodule connection $(E,\nabla_E, \sigma_E)$ as above, given by $\nabla_E(x) = 1 \otimes \extd{x}$ and $\sigma_E(\extd{g} \otimes x) = 1 \otimes \extd{g}.x$, gives a differential fibre bundle.
The fibre calculi are $N_q \cong (W^\perp)^{\wedge q}.\mathbb{C}X$, on which a differential $\extd : N_q \to N_{q+1}$ is given by
\begin{equation}
\extd(\xi.x)
= (-1)^{|\xi|} \xi \wedge \pi^{\perp}(\omega(x) \triangleleft x^{-1}).x
\end{equation}
for $\xi \in (W^{\perp})^{\wedge q}$ and $x \in X$.
The differential $\nabla_q : H^q(N) \to \Omega^1_{\mathbb{C}G} \otimes_{\mathbb{C}G} H^q(N)$ is given by
\begin{equation} \label{pruf}
\nabla_q([\xi.x])
= \pi(\omega(x) \triangleleft x^{-1}) \otimes [\xi.x].
\end{equation}
The fibre bundle $E$ gives rise to a spectral sequence converging to $H(\mathbb{C}X,E,\nabla_E) \cong H_{\rm{dR}}(\mathbb{C}X)$ with second page position $(p,q)$ given by $H^p(\mathbb{C}G, H^q(N), \nabla_{q})$ 
\end{proposition}
\begin{proof*}
\textbf{(1)} Firstly we show that $E$ is a differential fibre bundle.
The calculi $\Omega^p_B = \Omega^p_{\mathbb{C}G}$ are finitely freely generated for all $p \geq 0$ and therefore flat as modules, and the bimodule connection has zero curvature and is extendable.
Lastly we show that $g : \Omega^p_{\mathbb{C}G} \otimes_{\mathbb{C}G} M_{0,q} \to M_{p,q}$ given by $g(\xi \otimes [e \otimes \eta]) = [(\sigma_E \wedge \id)(\xi \otimes e \otimes \eta)] = [e \otimes \xi \wedge \eta]$ is an isomorphism.
Using $x \xi = x \xi x^{-1} x = (\xi \triangleleft x^{-1})x$ to move group elements to the right, and then $V = W \oplus W^{\perp}$;
\begin{align*}
M_{p,q} = \frac{\sigma_E(\Omega^p_{\mathbb{C}X} \otimes_{\mathbb{C}X} E) \wedge \Omega^q_{\mathbb{C}G}}{\sigma_E(\Omega^{p+1}_{\mathbb{C}X} \otimes_{\mathbb{C}X} E) \wedge \Omega^{q-1}_{\mathbb{C}G}}
= \frac{W^{\wedge p} \wedge V^{\wedge q}}{W^{\wedge p+1} \wedge V^{\wedge q-1}}.\mathbb{C}X
\cong W^{\wedge p} \otimes (W^{\perp})^{\wedge q}.\mathbb{C}X.
\end{align*}
The above isomorphism sends $[w_{i_1} \wedge \cdots \wedge w_{i_p} \wedge u_{j_1} \wedge \cdots \wedge u_{j_q}] \to w_{i_1} \wedge \cdots \wedge w_{i_p} \wedge u_{j_1} \wedge \cdots \wedge u_{j_q}$ where the $w_{i}$ are basis elements of $W$ and the $u_{i}$ are basis elements of $W^\perp$.
Thus $g$ sending $W^{\wedge p} \otimes (W^{\perp})^{\wedge q} \in \Omega^p_{\mathbb{C}G} \otimes_{\mathbb{C}X} M_{0,q}$ to $W^{\wedge p} \otimes (W^{\perp})^{\wedge q} \in M_{p,q}$ is an isomorphism.

\textbf{(2)}
The calculi on the fibres are $N_q \cong \frac{\Omega^q_{\mathbb{C}X}}{\Omega^1_{\mathbb{C}G} \wedge \Omega^{q-1}_{\mathbb{C}X}} \cong (W^\perp)^{\wedge q}.\mathbb{C}X$. From Lemma~\ref{lem111}
the differential $\extd : N_q \to N_{q+1}$ is given by $\extd(\xi.x) = (-1)^q [\xi \wedge \extd{x}]$ for $\xi \in (W^{\perp})^{\wedge q}$ and $x \in X$.
Using $\omega$ we rewrite this as
 $\extd(\xi.x) = (-1)^{|\xi|} \xi \wedge \pi^{\perp}(\omega(x) \triangleleft x^{-1}).x$.
 \newline\indent
If $\xi.x]\in H^q(N)$ then $\xi \wedge \pi^{\perp}(\omega(x) \triangleleft x^{-1}).x=0$, so
\[
\nabla_E^{(q)}( \xi.x )=(-1)^q\, \xi \wedge \pi(\omega(x) \triangleleft x^{-1}).x= \pi(\omega(x) \triangleleft x^{-1})  \wedge   \xi     .x
\]
By Proposition~\ref{prus} we have (\ref{pruf}). 
\end{proof*}


\subsection{Example: The group algebra $\mathbb{C} S_3$}

\begin{example}
Let $X = S_3$, denote transpositions as $u = (12)$ and $v = (23)$, and then define a subgroup $G = \{e, u\} \subset S_3$.
An example of a right representation of $X$ is given by $V = \mathbb{C}^2$ with right action $(v_1,v_2) \triangleleft x = (v_1,v_2) \rho(x)$ for the unitary representation $\rho : S_3 \to End(V)$ given by $\rho(u) = \begin{psmallmatrix} 1 & 0 \\ 0 & -1 \end{psmallmatrix}$and $\rho(v) = \tfrac{1}{2} \begin{psmallmatrix}-1 & \sqrt{3} \\ \sqrt{3} & 1 \end{psmallmatrix}$.
To define a calculus on $ \mathbb{C}S_3$ (and therefore by restriction a calculus on $\mathbb{C}G$) we need a cocycle $\omega : S_3 \to \mathbb{C}^2$ satisfying $\omega(xy) = \omega(x) \rho(y) + \omega(y)$.
For the cocycle to be a well-defined linear map, we need to be able to apply $\omega$ to the three relations of $S_3$, which are $u^2=e$, $v^2=e$, and $uvu=vuv$.
If we write $\omega(v) = (a,b)$ and $\omega(u) = (c,d)$, we have the following.

\textbf{(1)} Recalling that $\omega(e)=0$, the relation $u^2 = e$ gives:
\begin{align*}
0 = \omega(u^2) = \omega(u) \rho(u) + \omega(u) = (c,d) \begin{psmallmatrix} 1 & 0 \\ 0 & -1 \end{psmallmatrix} + (c,d) = (c,-d) + (c,d) = (2c,0).
\end{align*}
Hence $c = 0$.
We can normalise to get $d=1$ so that $\omega(u) = (0,1)$.

\textbf{(2)} The relation $v^2 = e$ gives:
\begin{align*}
0 &= \omega(v^2) = \omega(v) \rho(v) + \omega(v) = \omega(v)(\rho(v) + I_2) = (a,b) \tfrac{1}{2} \begin{psmallmatrix} 1 & \sqrt{3} \\ \sqrt{3} & 3 \end{psmallmatrix} = \tfrac{1}{2} (a + \sqrt{3}b , \sqrt{3}a + 3b)
\end{align*}
Both equations arising from this give that $a = - \sqrt{3}b$.
We already normalised when defining $\omega(u)$, so we simply have $b$ as a free parameter, giving $\omega(v) = (-\sqrt{3}b,b)$.

\textbf{(3)} Finally we have the relation $uvu = vuv$, and a calculation shows that this provides no additional restrictions. 
\newline
(Note that in \cite{S3-calculus}, a different calculus on $S_3$ is obtained by using the same  $\rho$ but on the representation $V = M_2(\mathbb{C})$ instead of $V = \mathbb{C}^2$.)

This gives a 1-parameter family of 2D calculi on $\mathbb{C}X$, with $\Lambda^1$ generated by $e_u = \omega(u) = (0,1)$ and $e_v = \omega(v) = b(-\sqrt{3}, 1)$, where $b \in \mathbb{C}$ is a free parameter.
Take the calculus on $\mathbb{C}G$ to be the vector space $W$ generated by $e_u$, so $\Omega^0_{\mathbb{C}G} = \mathbb{C}\{1,u\}$ and $\Omega^1_{\mathbb{C}G} = \omega(u).\mathbb{C}\{1,u\}$.

\quad We calculate the de Rham cohomology.
As long as $b \neq \tfrac{1}{2}$, this $\omega$ doesn't send any elements of $X$ other than $e$ to zero, the calculus on $\mathbb{C}X$ is connected.  Hence
(by \cite{QRG} Prop.\ 4.25) and as $\extd V=0$,
 $H_{dR}(\mathbb{C}X) = \Lambda_{\mathbb{C}X}$, the left invariant forms in each dimension.
Thus $H_{dR}^0(\mathbb{C}X) \cong \mathbb{C}$ and $H_{dR}^1(\mathbb{C}X) \cong \mathbb{C} \oplus \mathbb{C}$, while $H_{dR}^2(\mathbb{C}X) \cong \mathbb{C}$ has basis $\omega(u) \wedge \omega(v)$.

\quad We calculate the Leray-Serre spectral sequence explicitly for this example, where $E = A = \mathbb{C}X$, $B = \mathbb{C}G$. As $\omega(u) = (0,1)$ is a basis of $W$, it follows that $(1,0)$ is a basis of $W^{\perp}$, and hence $\pi^{\perp}(\lambda,\mu) = (\lambda,0)$.

From the formula above that $N_q \cong (W^\perp)^{\wedge q}.\mathbb{C}X$, we have $N_0 \cong \mathbb{C}X$ and $N_1 \cong (1,0).\mathbb{C}X$ and
all other $N_q$ are zero.
The one non-trivial differential is therefore $\extd : N_0 \to N_1$, given by $\extd{x} = \pi^{\perp}(\omega(x) \triangleleft x^{-1}).x$.

\quad Explicit calculation shows the kernel of $\extd : N_0 \to N_1$ is two-dimensional with basis elements $e$ and $u$.
Its image is four-dimensional with basis elements $(1,0).v$, $(1,0).uv$, $(1,0).vu$, $(1,0).uvu$.
Hence $H^0(N)$ is two-dimensional with basis elements $[e]$ and $[u]$, while $H^1(N)$ is two-dimensional with basis $[(1,0).e]$ and $[(1,0).u]$.

\quad The differential $\nabla_0 : H^0(N) \to \Omega^1_B \otimes_B H^0(N)$ is given on basis elements by 
\newline $\nabla_0([e]) = \pi(\omega(e) \triangleleft e^{-1}) \otimes [e] = 0$ 
 \newline $\nabla_0([u]) = \pi(\omega(u) \triangleleft u^{-1}) \otimes [u] = (0,-1) \otimes [u]$.

The differential $\nabla_1 : H^1(N) \to \Omega^1_B \otimes_B H^1(N)$ is given on basis elements by 
\newline $\nabla_1([(1,0).e]) = \pi(\omega(e) \triangleleft e^{-1}) \otimes [(0,1).e] = 0$ 
\newline  $\nabla_1([(0,1).u]) = \pi(\omega(u) \triangleleft u^{-1}) \otimes [(0,1).u] = (0,-1) \otimes [(0,1).u]$.

Hence $\nabla_0$ has kernel spanned by $[e]$ and image spanned by $(0,1).[u]$, while $\nabla_1$ has kernel spanned by $[(0,1).e]$ and image spanned by $(0,1) \otimes [(0,1).u]$.

As $\Omega^p_{\mathbb{C}X} = 0$ for $p \geq 2$ and $H^q(N) = 0$ for $q \geq 2$, the sequences for the cohomology are
\[
\begin{tikzcd}
0 \arrow{r}
&H^0(N) \arrow{r}{\nabla_0}
& \Omega^1_{\mathbb{C}G} \otimes_{\mathbb{C}G} H^0(N) \arrow{r}
& 0
\end{tikzcd}
\]
\[
\begin{tikzcd}
0 \arrow{r}
&H^1(N) \arrow{r}{\nabla_1}
& \Omega^1_{\mathbb{C}G} \otimes_{\mathbb{C}G} H^1(N) \arrow{r}
& 0
\end{tikzcd}
\]

$H^0(B,H^0(N),\nabla_0)$  is $\frac{\ker(\nabla_0)}{\mathrm{im}\, (0)} \cong \langle [e] \rangle_{\text{span}} \cong \mathbb{C}$.

$H^1(B,H^0(N),\nabla_0)$ is  $\frac{\Omega^1_{\mathbb{C}G} \otimes H^1(N)}{\mathrm{im}\, (\nabla_0)} \cong \langle (0,1) \otimes [e] \rangle_{\text{span}} \cong \mathbb{C}$.

$H^0(B,H^1(N),\nabla_1)$ is $\frac{\ker(\nabla_1)}{\mathrm{im}\, (0)} \cong \langle [(0,1).e] \rangle_{\text{span}} \cong \mathbb{C}$.

$H^1(B,H^1(N),\nabla_1)$ is $\frac{\Omega^1_{\mathbb{C}G} \otimes_{\mathbb{C}G}}{\mathrm{im}\,  \nabla_1} \cong \langle (1,0) \otimes [(1,0).e] \rangle_{\text{span}} \cong \mathbb{C}$.

\quad Page 2 of the Leray-Serre spectral sequence has entries $E^{p,q}_2 = H^p(\mathbb{C}G, H^q(N), \nabla_{q})$, with $E^{0,0}_2, E^{0,1}_2, E^{1,0}_2, E^{1,1}_2$ as its nonvanishing entries.
This is stable already, and hence the nontrivial cohomology groups are the following direct sums along diagonals.
\begin{align*}
&H^0(\mathbb{C}S_3,E,\nabla_E) \cong H^0(B, H^0(N), \nabla_{0}) \cong \mathbb{C} \\
&H^1(\mathbb{C}S_3,E,\nabla_E) \cong H^1(B, H^0(N), \nabla_{0}) \oplus H^0(B, H^1(N), \nabla_{1}) \cong \mathbb{C} \oplus \mathbb{C} \\
&H^2(\mathbb{C}S_3,E,\nabla_E) \cong H^1(B, H^1(N), \nabla_{1}) \cong \mathbb{C}
\end{align*}
This is the same as the de Rham cohomology $H_{dR}(\mathbb{C}X)$ that we calculated earlier.
\end{example}

\section{Example: Matrices} \label{fibre-bundles-example-matrices}

\begin{example} \label{algebra-index-M2-of-C} $\mathbf{M_2(\mathbb{C})}$
(\cite{QRG} Example 1.8, \cite{BeggsMajid17Matrices}.)
The algebra $M_2(\mathbb{C})$ of 2x2 complex-valued matrices has basis elements $E_{11}, E_{12}, E_{21}, E_{22}$ consisting of matrices with a 1 in the specified entry and all other entries zero.
It has an inner calculus given by $\extd{b} = [\theta', b] = \theta' b - b \theta'$ where
 $\theta' = E_{12}s' + E_{21}t'$ and $s'$ and $t'$ are central generators.
 The maximal prolongation calculus has the relation $s' \wedge t' = t' \wedge s'$.
\end{example}

We extend this idea to $M_3(\mathbb{C})$, giving it an inner calculus by $\theta = E_{12}s + E_{21}t + E_{33}u$ for central elements $s,t,u$.
The differential $\extd: M_3(\mathbb{C}) \to \Omega^1_{M_3(\mathbb{C})}$ is then given by $\extd{a} = [\theta, a] = [E_{12},a]s + [E_{21},a]t + [E_{33},a]u$, which on a general matrix in $M_3(\mathbb{C})$ is the following.
\begin{equation}
\extd \begin{psmallmatrix} a & b & c \\ d & e & f \\ g & h & i \end{psmallmatrix} =
\begin{psmallmatrix} d & e-a & f \\ 0 & -d & 0 \\ 0 & -g & 0 \end{psmallmatrix} s +
\begin{psmallmatrix} -b & 0 & 0 \\ a-e & b & c \\ -h & 0 & 0 \end{psmallmatrix} t +
\begin{psmallmatrix} 0 & 0 & -c \\ 0 & 0 & -f \\ g & h & 0 \end{psmallmatrix} u
\end{equation}
From this we can see that $\extd{E_{33}} = 0$, which means the calculus is not connected, since a connected calculus needs $\ker{\extd} = \mathbb{C}.I_3$.
For a higher order inner calculus, the differential is given by $\extd{\xi} = \theta \wedge \xi - (-1)^{|\xi|} \xi \wedge \theta$ for the inner element $\theta$.
For example, since $|u| = 1$, we have $\extd{u} = \theta \wedge u + u \wedge \theta$.

\begin{proposition}
Equipping $M_3(\mathbb{C})$ with higher order inner calculus for the inner element $\theta = E_{12}s + E_{21}t + E_{33}u$ necessitates that $s \wedge t = t \wedge s = u \wedge u$.
\end{proposition}

\begin{proof*}
As the calculus is inner, the differential is given by $\extd{a} = \theta a - a \theta$.
If we apply the differential twice to an element $a \in M_3(\mathbb{C})$, we get $\extd^2{a} = \theta \wedge (\theta a - a \theta) - (\theta a - a \theta) \wedge \theta = \theta \wedge \theta a - \theta \wedge a \theta + \theta \wedge a \theta - a \theta \wedge \theta = \theta \wedge \theta a - a \theta \wedge \theta = [\theta \wedge \theta , a]$.
For $\extd$ to be well-defined as a differential we need $\extd^2{a}$ to vanish, so $\theta \wedge \theta$ needs to be central so that its commutator with anything vanishes.
We calculate $\theta \wedge \theta = (E_{12}s + E_{21}t + E_{33}u) \wedge (E_{12}s + E_{21}t + E_{33}u) = E_{11}s \wedge t + E_{22} t \wedge s + E_{33} u \wedge u$.
The only central elements of $M_3(\mathbb{C})$ are multiples of $I_3$, and hence for $\theta \wedge \theta$ to be central we require $s \wedge t = t \wedge s = u \wedge u$.
\qed
\end{proof*}

Although the additional assumptions that $u \wedge t = t \wedge u$ and $u \wedge s = s \wedge u$ are not mandatory, we make these as well so that all the generators of the calculi commute.
Based on a private communication \cite{matrices-private-communication}, these extra assumptions bring the growth of the calculi down from exponential to polynomial.
With these additional assumptions, the derivatives of the calculi's basis elements are $\extd{s} = 2 s \wedge \theta$, $\extd{t} = 2 t \wedge \theta$ and $\extd{u} = 2 u \wedge \theta$.

For $A = M_3(\mathbb{C})$ and $B = M_2(\mathbb{C})$, $E = M_{2,3}(\mathbb{C})$ is a $B$-$A$ bimodule with actions given by matrix multiplication.

\begin{proposition} \label{pfus4}
Suppose $A$ has inner calculus as above given by inner element $\theta = E_{12}s + E_{21}t + E_{33}u$.
Then there is a unique right connection $\nabla_E : E \to E \otimes_A \Omega^1_A$ satisfying $\nabla_E(e_0) = 0$ for $e_0 = \begin{psmallmatrix} 2 & 0 & 0 \\ 0 & 2 & 0 \end{psmallmatrix}$. This is of the form 
$\nabla_E(e_0a) = e_0 \otimes \extd{a}$ and has zero curvature.
It is also an extendable bimodule connection by the bimodule map $\sigma_E : \Omega^1_B \otimes_B E \to E \otimes_A \Omega^1_A$ given by $\sigma_E(\extd \begin{psmallmatrix} a & b \\ c & d \end{psmallmatrix} \otimes e_0) = e_0 \otimes \extd \begin{psmallmatrix} a & b & 0 \\ c & d & 0 \\ 0 & 0 & 0 \end{psmallmatrix}$, which satisfies $\sigma_E(s' \otimes e_0) = e_0 \otimes s$ and $\sigma_E(t' \otimes e_0) = e_0 \otimes t$.
\end{proposition}

\begin{proof*}
\textbf{(1)} First we show well-definedness of $\nabla_E$.
Observing that
$e_0.\begin{psmallmatrix}
0 & 0 & 0 \\
0 & 0 & 0 \\
g & h & i
\end{psmallmatrix}
= \begin{psmallmatrix}
0 & 0 & 0 \\
0 & 0 & 0
\end{psmallmatrix} 
\in E$, the image of this under the linear map $\nabla_E$ must be zero, meaning that the differential must satisfy
$e_0 \otimes \extd \begin{psmallmatrix}
0 & 0 & 0 \\
0 & 0 & 0 \\
g & h & i
\end{psmallmatrix} = 0 \in E \otimes_A \Omega^1_A$.
We calculate using the differential above that
$e_0 \otimes \extd E_{3i} = e_0 \otimes \big(0 + 0 + [E_{33}, E_{3i}]u \big) = 
\begin{psmallmatrix}
2 & 0 & 0 \\
0 & 2 & 0
\end{psmallmatrix}
\otimes (E_{3i} - \delta_{i,3} E_{33})u = 
\begin{psmallmatrix}
2 & 0 & 0 \\
0 & 2 & 0
\end{psmallmatrix} (E_{3i} - \delta_{i,3} E_{33})
\otimes u = 0$, seeing as nonzero entries of $(E_{3i} - \delta_{i,3} E_{33})$ can only lie in the third row, and thus $\nabla_E$ is well-defined.

\textbf{(2)} Secondly, we calculate $\nabla_E$.
Every element of $E$ is of the form $e_0.a$, since $e_0.M_3(\mathbb{C}) = M_{2,3}(\mathbb{C}) = E$.
Therefore, using the Leibniz rule and the assumption $\nabla_E(e_0) = 0$, we calculate the connection as $\nabla_E(e_0a) = \nabla_E(e_0).a + e_0 \otimes \extd{a} = e_0 \otimes \extd{a}$.

\textbf{(3)} Thirdly, the map $\sigma_E$ satisfies $\sigma_E(\extd{b} \otimes e_0) = \nabla_E(b e_0) - b \nabla_E(e_0)$.
But $\nabla_E(e_0) = 0$, so $\sigma_E(\extd \begin{psmallmatrix} a & b \\ c & d \end{psmallmatrix} \otimes e_0) = \nabla_E(\begin{psmallmatrix} a & b \\ c & d \end{psmallmatrix} e_0) = \nabla_E(e_0 \begin{psmallmatrix} a & b & 0 \\ c & d & 0 \\ 0 & 0 & 0 \end{psmallmatrix}) = e_0 \otimes \extd \begin{psmallmatrix} a & b & 0 \\ c & d & 0 \\ 0 & 0 & 0 \end{psmallmatrix}$ as required.

\textbf{(4)} Next, we show $\sigma_E(s' \otimes e_0) = e_0 \otimes s$ and $\sigma_E(t' \otimes e_0) = e_0 \otimes t$.
In the calculus on $B$, we have $\extd{E_{21}} = [E_{12}, E_{21}]s' + [E_{21},E_{21}]t' = \begin{psmallmatrix} 1 & 0 \\ 0 & -1 \end{psmallmatrix}s'$, and likewise on the calculus on $A$.
Therefore, using the fact that $\sigma_E$ is a bimodule map and that $s'$ is central and also the formula above for $\sigma_E$,
\begin{align*}
&\begin{psmallmatrix} 1 & 0 \\ 0 & -1 \end{psmallmatrix} \sigma_E(s' \otimes e_0)
= \sigma_E(s' \begin{psmallmatrix} 1 & 0 \\ 0 & -1 \end{psmallmatrix} \otimes e_0)
= \sigma_E(\extd{E_{21}} \otimes e_0) \\
&= {e_0 \otimes \extd \begin{psmallmatrix} 0 & 0 & 0 \\ 1 & 0 & 0 \\ 0 & 0 & 0 \end{psmallmatrix}}
= e_0 \otimes \begin{psmallmatrix} 1 & 0 & 0 \\ 0 & -1 & 0 \\ 0 & 0 & 0 \end{psmallmatrix} s
= \begin{psmallmatrix} 1 & 0 \\ 0 & -1 \end{psmallmatrix} e_0 \otimes s.
\end{align*}
However, as $\begin{psmallmatrix} 1 & 0 \\ 0 & -1 \end{psmallmatrix}$ is invertible, this implies $\sigma_E(s' \otimes e_0) = e_0 \otimes s$.
The result $\sigma_E(t' \otimes e_0) = e_0 \otimes t$ follows similarly by considering $\extd{E_{12}} = [E_{12}, E_{12}]s' + [E_{21},E_{12}]t' = \begin{psmallmatrix} -1 & 0 \\ 0 & 1 \end{psmallmatrix}t'$.

\textbf{(5)} 
Since $B = M_2(\mathbb{C})$ is equipped with maximal prolongation calculus, Corollary 5.3 of \cite{extendability-for-bimodules} says that every zero-curvature bimodule connection is extendable.
\qed
\end{proof*}

\begin{proposition}
Suppose $B = M_2(\mathbb{C})$ and $A = M_3(\mathbb{C})$ are equipped with the above calculi.
Then the $B$-$A$ bimodule $E = M_{2,3}(\mathbb{C})$ with the bimodule connection $(\nabla_E, \sigma_E)$ from Proposition~\ref{pfus4} gives a differential fibre bundle, and thus a spectral sequence converging to $H(A,E,\nabla_E) = H(M_3(\mathbb{C}), M_{2,3}(\mathbb{C}), \nabla_E)$.
\end{proposition}

\begin{proof*}
For all $p \geq 0$ the calculi $\Omega^p_B = \Omega^p_{M_2(\mathbb{C})}$ are freely finitely generated and hence flat as modules.
The bimodule connection $(\nabla_E, \sigma_E)$ satisfies the requirements of having zero curvature and being extendable.
The last property we need to show is that the map $g : \Omega^p_B \otimes_B M_{0,q} \to M_{p,q}$ given by $g(\xi \otimes [e \otimes \eta]) = [(\sigma_E \wedge \id)(\xi \otimes e \otimes \eta)]$ is an isomorphism.

Since $E = e_0.M_3(\mathbb{C})$, the forms on $M_3(\mathbb{C})$ of degree $p$ in the fibre and $q$ in the base are 
\begin{align*}
M_{p,q} = \frac{\sigma_E(\Omega^p_B \otimes_B E) \wedge \Omega^q_A}{\sigma_E(\Omega^{p+1}_B \otimes_B E) \wedge \Omega^{q-1}_A}
\cong \frac{\sigma_E(\Omega^p_{M_2(\mathbb{C})} \otimes_{M_2(\mathbb{C})} e_0) \wedge \Omega^q_{M_3(\mathbb{C})}}{\sigma_E(\Omega^{p+1}_{M_2(\mathbb{C})} \otimes_{M_2(\mathbb{C})} e_0) \wedge \Omega^{q-1}_{M_3(\mathbb{C})}}.
\end{align*}

Everything in the numerator is of the form $(s \text{ or } t)^{\wedge(p+q-k)} \wedge u^{\wedge k} .M_3(\mathbb{C})$ for some $0 \leq k \leq q$, while everything in the denominator is of the form $(s \text{ or } t)^{\wedge(p+q-k+1)} \wedge u^{\wedge (k-1)} .M_3(\mathbb{C})$ for $0 \leq k \leq q$.
Since $u \wedge u = s \wedge t$, it follows that if an element of the numerator has $k \geq 2$ then it lies in the denominator.
But if an element of the numerator has $k < q$ then it has to lie in the denominator.
Therefore $M_{p,q} = 0$ for $q \geq 2$, and (omitting to write the equivalence classes) a basis of $M_{p,0}$ is given by $e_0 \otimes s^{\wedge r} \wedge t^{\wedge(p-r)}$ for some $0 \leq r \leq p$, while a basis of $M_{p,1}$ is given by $e_0 \otimes s^{\wedge r} \wedge t^{\wedge(p-r)} \wedge u$.

In the case $q=0$, the map $g$ is given on basis elements as
\begin{align*} 
(s')^{\wedge r} \wedge (t')^{\wedge(p-r)} \otimes e_0 \longmapsto e_0 \otimes s^{\wedge r} \wedge t^{\wedge(p-r)}.
\end{align*}
The map $g$ here is an isomorphism, since it just re-arranges the order of the tensor product and and re-labels $s'$ and $t'$, which introduces no new relations.

Similarly in the case $q=1$, the map $g$ is given on basis elements as
\begin{align*}
(s')^{\wedge r} \wedge (t')^{\wedge(p-r)} \otimes e_0 \otimes u \longmapsto e_0 \otimes s^{\wedge r} \wedge t^{\wedge(p-r)} \wedge u,
\end{align*}
which is an isomorphism.
\qed
\end{proof*}

\begin{proposition}
The nonzero cohomology groups of $A$ with coefficients in $E$ can be calculated via the Leray-Serre spectral sequence as $H^0(A,E,\nabla_E) \cong  \mathbb{C}$, $H^1(A,E,\nabla_E) \cong \mathbb{C}^{6}$, $H^2(A,E,\nabla_E) \cong \mathbb{C}^{5}$.
\end{proposition}

\begin{proof*}
\textbf{(1)} Find $N_i$ : The space $N_0$ is isomorphic to $e_0.A$, which has six-dimensional vector space basis $e_0.E_{ij}$ for $1 \leq i \leq 2$ and $1 \leq j \leq 3$ (i.e.\ excluding the bottom row).
The space $N_1$ is isomorphic to $e_0.A \otimes u$, which has six-dimensional vector space basis $e_0.E_{ij} \otimes u$ for $1 \leq i \leq 2$ and $1 \leq j \leq 3$.
Since we showed earlier that $M_{p,q} = 0$ for $q \geq 2$, this means $N_i = M_{0,i} = 0$ for $i \geq 2$.

\textbf{(2)} The differential $\extd : N_0 \to N_1$ is given by $\extd([e_0.E_{ij}]) = [\nabla_E(e_0.E_{ij})] = [e_0 \otimes \extd{E_{ij}}] = [e_0.[E_{33}, E_{ij}] \otimes u]$.
The kernel has four-dimensional basis $[e_0.E_{11}]$, $[e_0.E_{12}]$, $[e_0.E_{21}]$, $[e_0.E_{22}]$.
The image has two-dimensional basis $[e_0.E_{13} \otimes u]$ and $[e_0.E_{23} \otimes u]$.

\textbf{(3)} Consequently $H^0(N)$ is four-dimensional with basis $[[e_0.E_{11}]]$, $[[e_0.E_{12}]]$, $[[e_0.E_{21}]]$, $[[e_0.E_{22}]]$.
Also $H^1(N)$ is four-dimensional with basis $[[e_0.E_{11} \otimes u]]$, $[[e_0.E_{12} \otimes u]]$, $[[e_0.E_{21} \otimes u]]$, $[[e_0.E_{22} \otimes u]]$.

\textbf{(4)} Calculate $\nabla_0 : H^0(N) \to \Omega^1_B \otimes_B H^0(N)$ on the basis of $H^0(N)$.
For $1 \leq i,j \leq 2$ we have $\nabla_0([[e_0.E_{ij}]]) = g^{-1}([[e_0 \otimes \extd{E_{ij}}]])$.
Then
\begin{align*}
& \nabla_0(e_0.E_{12}) = g^{-1}([[e_0 \otimes (E_{22} - E_{11})t]]) = t' \otimes [[e_0.(E_{22} - E_{11})]], \\
&\nabla_0(e_0.E_{21}) = s' \otimes [[e_0.(E_{11} - E_{22})]], \\
& \nabla_0(e_0.E_{11}) = -s' \otimes [[e_0.E_{12}]] + t' \otimes [[e_0.E_{21}]] = - \nabla_0(e_0.E_{22}).
\end{align*}
Hence $\nabla_0$ has one-dimensional kernel with basis $[[e_0.(E_{11} + E_{22})]]$, and three-dimensional image with basis elements $t' \otimes [[e_0.(E_{22} - E_{11})]]$, $s' \otimes [[e_0.(E_{11} - E_{22})]]$, $t' \otimes [[e_0.E_{21}]] - s' \otimes [[e_0.E_{12}]]$.

\textbf{(5)} Calculate $\nabla_1 : H^1(N) \to \Omega^1_B \otimes_B H^1(N)$ on the basis of $H^1(N)$.
For $1 \leq i,j \leq 2$,
\begin{align*}
\nabla_E^{[1]}(e_0.E_{ij} \otimes u) &= \nabla_E(e_0.E_{ij}) \wedge u + e_0.E_{ij} \otimes \extd{u} = e_0 \otimes [\theta, E_{ij}] \wedge u + 2 e_0.E_{ij} \otimes \theta \wedge u \\
&= e_0 \otimes \Big( (E_{12} E_{ij} + E_{ij} E_{12}) s + (E_{21} E_{ij} + E_{ij} E_{21}) t \Big) \wedge u \\
&= \sigma_E(s' \otimes e_0.(E_{12} E_{ij} + E_{ij} E_{12})) \wedge u + \sigma_E(t' \otimes e_0.(E_{21} E_{ij} + E_{ij} E_{21})) \wedge u.
\end{align*}

Consequently, 
\begin{align*}
\nabla_1([[e_0.E_{ij} \otimes u]]) = s' \otimes [[e_0.(E_{12} E_{ij} + E_{ij} E_{12}) \otimes u]] + t' \otimes [[e_0.(E_{21} E_{ij} + E_{ij} E_{21}) \otimes u]].
\end{align*}

Using this, 
$\nabla_1([[e_0.E_{12} \otimes u]]) = t' \otimes [[e_0 \otimes u]]$
and
$\nabla_1([[e_0.E_{21} \otimes u]]) = s' \otimes [[e_0 \otimes u]]$
and
$\nabla_1([[e_0.E_{11} \otimes u]]) = s' \otimes [[e_0.E_{12} \otimes u]] + t' \otimes [[e_0.E_{21} \otimes u]] = \nabla_1([[e_0.E_{22} \otimes u]])$.

Hence the kernel of $\nabla_1$ has one-dimensional basis $[[e_0.(E_{11} - E_{22}) \otimes u]]$, while the image has three-dimensional basis $t' \otimes [[e_0 \otimes u]]$, $s' \otimes [[e_0 \otimes u]]$, $s' \otimes [[e_0.E_{12} \otimes u]] + t' \otimes [[e_0.E_{21} \otimes u]]$.

\textbf{(6)} Next we work out the quotients for cohomology.

\textbf{(i)} Firstly, $H^0(B,H^0(N),\nabla_0) \cong \frac{\ker(\nabla_0)}{\mathrm{im}\, (0)} \cong \mathbb{C}$.

\textbf{(ii)} Secondly, $H^1(B,H^0(N),\nabla_0) \cong \frac{\Omega^1_{B} \otimes H^0(N)}{\mathrm{im}\, (\nabla_0)}$.
As $\Omega^1_B$ is a free module with two basis elements and $H^0(N)$ is four-dimensional, the vector space $\Omega^1_B \otimes_B H^0(N)$ is eight-dimensional.
The quotient is therefore five dimensional, and an example of a basis of $\frac{\Omega^1_{B} \otimes H^0(N)}{\rm{im}(\nabla_0)}$ is given by
$[s' \otimes [[e_0.E_{11}]]]$,
$[s' \otimes [[e_0.E_{12}]]]$,
$[s' \otimes [[e_0.E_{21}]]]$,
$[t' \otimes [[e_0.E_{11}]]]$,
$[t' \otimes [[e_0.E_{12}]]]$.
Hence $H^1(B,H^0(N),\nabla_0) \cong \mathbb{C}^{5}$.

\textbf{(iii)} Thirdly, $H^0(B,H^1(N),\nabla_1) \cong \frac{\ker(\nabla_1)}{\mathrm{im}\, (0)} \cong \mathbb{C}$.

\textbf{(iv)} Lastly, $H^1(B,H^1(N),\nabla_1) \cong \frac{\Omega^1_{B} \otimes_{B} H^1(N)}{\mathrm{im}\, (\nabla_1)}.$
As $\Omega^1_B$ is a free module with two basis elements and $H^1(N)$ is four-dimensional, the vector space $\Omega^1_B \otimes_B H^1(N)$ is eight-dimensional.
Taking the quotient by the three-dimensional $\rm{im}(\nabla_1)$ gives a five-dimensional vector space.
Hence $H^1(B,H^1(N),\nabla_1) \cong \mathbb{C}^{5}$.

\textbf{(7)} Page 2 of the Leray-Serre spectral sequence has entries $E^{p,q}_2 = H^p(B, H^q(N), \nabla_{q})$, with $E^{0,0}_2, E^{0,1}_2, E^{1,0}_2, E^{1,1}_2$ as its nonvanishing entries.
This is stable already, and hence the nontrivial cohomology groups are the following direct sums along diagonals.
\begin{align*}
&H^0(A,E,\nabla_E) \cong H^0(B, H^0(N), \nabla_{0}) \cong \mathbb{C} \\
&H^1(A,E,\nabla_E) \cong H^1(B, H^0(N), \nabla_{0}) \oplus H^0(B, H^1(N), \nabla_{1}) \cong \mathbb{C}^{5} \oplus \mathbb{C} \cong \mathbb{C}^{6} \\
&H^2(A,E,\nabla_E) \cong H^1(B, H^1(N), \nabla_{1}) \cong \mathbb{C}^{5}. \qquad\qquad \qquad\qquad \square
\end{align*}

\end{proof*}

The bimodule $E$ has inner product $\langle , \rangle : \overline{E} \otimes_B E \to A$ given by $\langle \overline{x} , y \rangle = x^*y$, where $*$ is the conjugate transpose map.
As matrix algebras are C*-algebras, the KSGNS construction says that the map $\phi : B \to A$ given by $\phi(b) = \langle \overline{e_0} , b e_0 \rangle$ is completely positive.
For $e_0 = \begin{psmallmatrix} 2 & 0 & 0 \\ 0 & 2 & 0 \end{psmallmatrix}$ then $\phi$ is not an algebra map, as $\phi(I_2) $ is not a multiple of $ I_3$, and algebra maps have to send the identity to the identity.

Moreover, $\nabla_E(e_0) = 0$, so for $\phi$ to be a cochain map we just need metric preservation, which holds because of the following.
Recall that for the right connection $\nabla_E$ on $E$, we have $\nabla_E(e_0a) = e_0 \otimes \extd{a}$, which gives a corresponding left connection $\nabla_{\overline{E}}$ on $\overline{E}$ given by $\nabla_{\overline{E}}(\overline{e_0a}) = \extd{a}^* \otimes \overline{e_0}$, and that the inner product on $E$ is given by $\langle \overline{x} , y \rangle = x^{-1}y$.
Then:
\begin{align*}
&\extd{a_1}^* \langle \overline{e_0} , e_0 a_2 \rangle + \langle \overline{e_0 a_1} , e_0 \rangle \extd{a_2} = \extd{a_1}^*.e_0^* e_0 a_2 + a_1^* e_0^* e_0 \extd{a_2} = 4 \extd{a_1}^*.a_2 + 4a_1^* \extd{a_2} \\
&= 4 \extd{a_1^* a_2} = \extd{a_1^* e_0^* e_0 a_2} = \extd \langle \overline{e_0 a_1} , e_0 a_2 \rangle.
\end{align*}
Thus by Proposition 4.86 of \cite{QRG}, $\phi$ is a completely positive cochain map, but not an algebra map.

\section{Example: Noncommutative torus fibred over the circle} \label{fibre-bundles-example-circle-in-torus}

We look at an infinite dimensional example, with total algebra the quantum torus $A = \mathbb{C}_{\theta}[\mathbb{T}^2]$ and base  algebra the circle $B = \mathbb{C}[S^1]$ with its classical calculus. Here 
 $\mathbb{C}[S^1] = \mathbb{C}[t,t^{-1}]$ is the algebra of polynomials in $t$ and $t^{-1}$
 with $*$-structure given by $t^* = t^{-1}$.

\begin{example} \label{algebra-index-noncommutative-torus} $\mathbf{\mathbb{C}_{\theta}[\mathbb{T}^2]}$
(\cite{QRG} Example 1.36. Calculi on the noncommutative torus were first studied by Connes and Rieffel in \cite{ConnesRieffel87}.)
The noncommutative torus $\mathbb{C}_{\theta}[\mathbb{T}^2]$ has generators $u,v$ with the relation $vu = e^{i\theta}uv$ for a real parameter $\theta$.
It has star structure $u^* = u^{-1}$ and $v^* = v^{-1}$.
It has a calculus $\Omega^1 = \mathbb{C}_{\theta}[\mathbb{T}^2].\{ \extd{u}, \extd{v} \}$, with left action given by multiplication and right action given by
\begin{align*}
\extd{u}.u = u.\extd{u}, \quad
\extd{v}.v = v.\extd{v}, \quad
\extd{v}.u = e^{i\theta}u.\extd{v}, \quad
\extd{u}.v = e^{-i\theta}v.\extd{u}.
\end{align*}
We take maximal prolongation calculi for the higher calculi on $A$, giving relations $\extd{u} \wedge \extd{u} = 0 = \extd{v} \wedge \extd{v}$ and $\extd{v} \wedge \extd{u} = - e^{i \theta} \extd{u} \wedge \extd{v}$.
These relations imply $\Omega^2_A = \extd{u} \wedge \extd{v}.A$ and $\Omega^q_A = 0$ for $q \geq 3$, with every nonzero element of $\Omega^2_A$ a multiple of $\extd{u} \wedge \extd{v}$.
\end{example}

We take a $B$-$A$ bimodule $E$ given by $E = \mathbb{C}_{\theta}[\mathbb{T}^2] \oplus \mathbb{C}_{\theta}[\mathbb{T}^2]$, with left $B$-action and right $A$-action given respectively by
\begin{align*}
t \triangleright(f \oplus g) = uf \oplus vg, \quad
(f \oplus g) \triangleleft k = fk \oplus g k\ .
\end{align*}

\begin{proposition} \label{TorusConnection}
If $B = \mathbb{C}_q[S^1]$ is equipped with classical calculus ($q=1$), then there is a zero curvature right connection $(\nabla_E, \sigma_E)$ on $E$, for  $\nabla_E : E \to E \otimes_A \Omega^1_A$ and  $\sigma_E : \Omega^1_B \otimes_B E \to E \otimes_A \Omega^1_A$
given by
\begin{align*}
\nabla_E(f \oplus 0) &= (1 \oplus 0) \otimes \extd{f}\ , \quad
\nabla_E(0 \oplus g) = (0 \oplus 1) \otimes \extd{g}\ ,\cr
\sigma_E(\extd{t} \otimes (f \oplus 0)) &= (1 \oplus 0) \otimes \extd{u}.f\ , \quad
\sigma_E(\extd{t} \otimes (0 \oplus g)) = (0 \oplus 1) \otimes \extd{v}.g\ .
\end{align*}
Since higher calculi on the circle are zero, this connection is automatically extendable.
\end{proposition}

\begin{proof*}
\textbf{(1)} Firstly, we calculate the formula for $\sigma_E$ on the generators.
\begin{align*}
\sigma_E&(\extd{t} \otimes (f \oplus 0)) = \nabla_E(t.(f \oplus 0)) - t.\nabla_E(f \oplus 0) \\
&= \nabla_E(uf \oplus 0) - t.(1 \oplus 0) \otimes \extd{f} = (1 \oplus 0) \otimes \extd(uf) - (u \oplus 0) \otimes \extd{f} \\
&= (1 \oplus 0) \otimes \extd{u}.f + (1 \oplus 0) \otimes u.\extd{f} - (u \oplus 0) \otimes \extd{f} = (1 \oplus 0) \otimes \extd{u}.f
\end{align*}
Here we use the standard formula for a bimodule connection $\sigma(\extd{a} \otimes e) = \nabla(a.e) - a.\nabla(e)$.
\begin{align*}
\sigma_E&(\extd{t} \otimes (0 \oplus g)) = \nabla_E(t.(0 \oplus g)) - t.\nabla_E(0 \oplus g) \\
&= \nabla_E(0 \oplus vg) - (0 \oplus v) \otimes \extd{g} = (0 \oplus 1) \otimes \extd(vg) - (0 \oplus v) \otimes \extd{g} \\
&= (0 \oplus 1) \otimes v.\extd{g} + (0 \oplus 1) \extd{v}.g - (0 \oplus v) \otimes \extd{g} = (0 \oplus 1) \otimes \extd{v}.g
\end{align*}

\textbf{(2)} Next we show this is a bimodule map.
The calculations in part (1) show that $\sigma_E$ is a right module map, so we just need to show it is a left module map.
Using $t.\extd{t} =  \extd{t}.t$ and $\extd{u}.u = u.\extd{u}$, we calculate:
\begin{align*}
\sigma_E&(t.\extd{t} \otimes (1 \oplus 0)) =  \sigma_E(\extd{t} \otimes t.(1 \oplus 0)) 
= \sigma_E(\extd{t} \otimes (u \oplus 0)) =  (1 \oplus 0) \otimes \extd{u}.u \\
&=  (1 \oplus 0) \otimes u.\extd{u} = (u \oplus 0) \otimes \extd{u} = t \triangleright (1 \oplus 0) \otimes \extd{u} \ . \qquad\square
\end{align*}
\end{proof*}

Having a zero curvature extendable bimodule connection, we show the last ingredient required for a bimodule fibre bundle.

\begin{proposition}
The map $g : \Omega^p_B \otimes_B M_{0,q} \to M_{p,q}$ given by
\begin{equation}
g(\xi \otimes [e \otimes \eta]) = [(\sigma_E \wedge \id)(\xi \otimes e \otimes \eta)]
\end{equation}
is an isomorphism, where forms of degree $p$ in the fibre and $q$ in the base are given by the formula:
\begin{equation}
M_{p,q} = \frac{\sigma_E(\Omega^p_B \otimes_B E) \wedge \Omega^q_A}{\sigma_E(\Omega^{p+1}_B \otimes_B E) \wedge \Omega^{q-1}_A}
\end{equation}
\end{proposition}

\begin{proof*}
We get, for $A= \mathbb{C}_{\theta}[\mathbb{T}^2]$
\begin{align*}
M_{0,0} &=N_0= E= \mathbb{C}_{\theta}[\mathbb{T}^2] \oplus \mathbb{C}_{\theta}[\mathbb{T}^2]\ ,\quad 
M_{1,0} = (1\oplus 0)\otimes \extd u\,  \mathbb{C}_{\theta}[\mathbb{T}^2] +
(0\oplus 1)\otimes \extd v\,  \mathbb{C}_{\theta}[\mathbb{T}^2] \ ,\cr
M_{0,1} &= N_1= \frac{ E\otimes \Omega^1_A}{(1\oplus 0)\otimes \extd u\,  \mathbb{C}_{\theta}[\mathbb{T}^2] +
(0\oplus 1)\otimes \extd v\,  \mathbb{C}_{\theta}[\mathbb{T}^2] }\ , \quad M_{1.1} = E\otimes \Omega^2_A\ .
\end{align*}
Then $g : \Omega^1_B \otimes_B M_{0,0} \to M_{1,0}$ gives
\begin{align*}
g(\extd{t} \otimes (1 \oplus 0)) = (1 \oplus 0) \otimes \extd{u}, \quad
g(\extd{t} \otimes (0 \oplus 1)) = (0 \oplus 1) \otimes \extd{v}
\end{align*}
is an isomorphism. Also $g : \Omega^1_B \otimes_B M_{0,1} \to M_{1,1}$ is an isomorphism, as
\begin{align*}
g(\extd{t} \otimes [(1 \oplus 0) \otimes \extd{v}.A] ) &= (1 \oplus 0) \otimes \extd{u} \wedge \extd{v}.A\ ,\cr
g(\extd{t} \otimes [(0 \oplus 1) \otimes \extd{u}.A ] ) &= (0 \oplus 1) \otimes \extd{v} \wedge \extd{u}.A
= (0 \oplus 1) \otimes \extd{u} \wedge \extd{v}.A \qquad\square
\end{align*}
\end{proof*}

Having shown that $g$ is an isomorphism, we have a bimodule differential fibre bundle, so there exists a Leray-Serre spectral sequence converging to the sheaf cohomology of the quantum torus with coefficients in the bimodule $E$.

\begin{proposition}
The sheaf cohomology given by the Leray-Serre spectral sequence for $E = \mathbb{C}_{\theta}[\mathbb{T}^2] \oplus \mathbb{C}_{\theta}[\mathbb{T}^2]$ and the connection $\nabla_E$ from Proposition \ref{TorusConnection} is:
\begin{align*}
H^0(\mathbb{C}_{\theta}[\mathbb{T}^2],E,\nabla_E) \cong \mathbb{C}^2, \quad
H^1(\mathbb{C}_{\theta}[\mathbb{T}^2],E,\nabla_E) \cong \mathbb{C}^{2} \oplus \mathbb{C}^2, \quad
H^2(\mathbb{C}_{\theta}[\mathbb{T}^2],E,\nabla_E) \cong \mathbb{C}^{2}. \qquad\qquad \qquad\qquad \square
\end{align*}
\end{proposition}

\begin{proof*}
We use the noncommutative partial derivative notation on $\mathbb{C}_{\theta}[\mathbb{T}^2],$
\[
\extd f = \extd u\, \partial_u f +  \extd v\, \partial_v f 
\]
and then $\extd :N_0\to N_1$ is
\[
\extd(f\oplus g) = [(1\oplus 0)\otimes \extd v\, \partial_v f ] + [(0\oplus 1)\otimes 
 \extd u\, \partial_u g ] \ .
\]
Then
\begin{align*}
H^0(N) &=\mathbb{C}[u,u^{-1}] \oplus \mathbb{C}[v,v^{-1}]\ ,\cr
H^1(N) &=
 [(1\oplus 0)\otimes \extd v\, v^{-1}\mathbb{C}[u,u^{-1}]]
 +
 [(0\oplus 1)\otimes  \extd u\,  u^{-1}\mathbb{C}[v,v^{-1}]]
\end{align*}
which are $\mathbb{C}[S^1]$ modules as required. Then
$\nabla_i:H^i(N) \to \Omega^1_B\otimes_B H^i(N)
$
are given by
\begin{align*}
\nabla_0&(f\oplus g) = \extd t \otimes (\partial_u f \oplus \partial_v g)\ ,\cr
\nabla_1&(  [(1\oplus 0)\otimes \extd v\, v^{-1}\,f  ] ) = - g^{-1}(
(1\oplus 0)\otimes \extd v\, v^{-1}\wedge \extd u\, \partial_u f) \cr
& =  g^{-1}(
(1\oplus 0)\otimes \extd u\wedge \extd v\, v^{-1}  \partial_u f) =
\extd t \otimes [ (1\oplus 0)\otimes \extd v\, v^{-1}  \partial_u f]    \ ,\cr
\nabla_1&(  [(0\oplus 1)\otimes  \extd u\,  u^{-1}  g] ) =
\extd t \otimes [ (0\oplus 1)\otimes \extd u\, u^{-1}  \partial_v g ] 
\end{align*}
for $f\in \mathbb{C}[u,u^{-1}]$ and $g\in \mathbb{C}[v,v^{-1}]$. Then \smallskip  \newline 
$H^0(B, H^0(N), \nabla_{0}) $ has basis
$1\oplus 0\ ,\ 0\oplus 1$\ ,
\newline 
$H^1(B, H^0(N), \nabla_{0}) $ has basis
$t^{-1}\extd t \otimes (1\oplus 0)\ ,\ t^{-1}\extd t \otimes (0\oplus 1)$\ ,
\newline 
$H^0(B, H^1(N), \nabla_{1}) $ has basis
$ [(0\oplus 1)\otimes  \extd u\,  u^{-1}  ]\ ,\ 
 [ (1\oplus 0)\otimes \extd v\, v^{-1}  ]  $\ ,
 \newline 
$H^1(B, H^1(N), \nabla_{1}) $ has basis
$\extd t \otimes [ (1\oplus 0)\otimes \extd v\, v^{-1} u^{-1}] 
=t^{-1} \extd t \otimes [ (1\oplus 0)\otimes \extd v\, v^{-1} ]   $\ ,\ 
$t^{-1}\extd t \otimes [ (0\oplus 1)\otimes \extd u\, u^{-1}  ] $\ .

The Leray-Serre spectral sequence is stable already from page 2, and hence the nontrivial cohomology groups are the following direct sums along the diagonals:
\begin{align*}
&H^0(A,E,\nabla_E) \cong H^0(B, H^0(N), \nabla_{0}) \cong \mathbb{C}^2 \\
&H^1(A,E,\nabla_E) \cong H^1(B, H^0(N), \nabla_{0}) \oplus H^0(B, H^1(N), \nabla_{1}) \cong \mathbb{C}^{2} \oplus \mathbb{C}^2 \\
&H^2(A,E,\nabla_E) \cong H^1(B, H^1(N), \nabla_{1}) \cong \mathbb{C}^{2}.
\end{align*}
\qed
\end{proof*}

\begin{remark}
In Example 1.36 of \cite{QRG}, the standard de Rham cohomology of the quantum torus is given as $H^0_{dR}(\mathbb{C}_{\theta}[\mathbb{T}^2]) \cong \mathbb{C}$, $H^1_{dR}(\mathbb{C}_{\theta}[\mathbb{T}^2]) \cong \mathbb{C} \oplus \mathbb{C}$, and $H^2_{dR}(\mathbb{C}_{\theta}[\mathbb{T}^2]) \cong \mathbb{C}$.
And so we can see that taking sheaf cohomology with coefficients in $E = \mathbb{C}_{\theta}[\mathbb{T}^2] \oplus \mathbb{C}_{\theta}[\mathbb{T}^2]$ gives the same result but with each instance of $\mathbb{C}$ replaced by $\mathbb{C}^2$.
\end{remark} 

Suppose we equip $E$ with inner product $\langle , \rangle : \overline{E} \otimes_B E \to A$ given by $\langle \overline{f_1 \oplus g_1} , f_2 \oplus g_2 \rangle = f_1^*f_2 + g_1^*g_2$.
Then $\nabla_E$ preserves this inner product and satisfies the conditions of Proposition~4.86 of \cite{QRG}. 
In the case $e_0 = 1 \oplus 1$, we have $\nabla_E(e_0) = 0$. Thus the 
 completely positive map $\phi : B \to A$ given by the KSGNS construction \cite{Lance} as $\phi(t^n) = \langle \overline{e_0} , t^n e_0 \rangle = u^n + v^n$ gives a cochain map but is not an algebra map. This qualifies as a fibre bundle under our definition which generalises the algebra map definition.


\begin{thebibliography}{99}
  \setlength{\parskip}{0pt}
  \setlength{\itemsep}{0pt plus 0.3ex}
  
  \small
  
  \bibitem{BrzezinskiSzymanski21}T. Brzezi\'{n}ski and W. Szyma\'{n}ski, An algebraic framework for noncommutative bundles with homogeneous fibres, Alg. Number Th. 15 (2021) 217-240, \url{https://doi.org/10.2140/ant.2021.15.217}

  
\bibitem{extendability-for-bimodules} G. Alhamzi and E.J. Beggs, Differentiating the State Evaluation Map from Matrices to Functions on Projective Space, Symmetry 2023, 15, 474. \url{https://doi.org/10.3390/sym15020474}
\bibitem{bbsheaf} E.J. Beggs and  T. Brzezi\'nski, The Serre spectral sequence of a noncommutative fibration for de Rham cohomology, Acta Math. 195 (2005) 155--196 \url{https://doi.org/10.1007/BF02588079}
\bibitem{QRG} E.J. Beggs and S. Majid, Quantum Riemannian Geometry, ISBN 978-3-030-30293-1, Springer, 2020 \url{https://doi.org/10.1007/978-3-030-30294-8}
\bibitem{BeggsMajid17Matrices} E.J. Beggs and S. Majid, Spectral triples from bimodule connections and Chern connections, J. Noncom. Geom. 11 (2017) 669–701 \url{https://doi.org/10.4171/jncg/11-2-7}
\bibitem{Bar-categories} E.J. Beggs and S. Majid, Bar categories and star operations, Algebras and Representation Theory, 12 (2009) 103-152 \url{https://doi.org/10.1007/s10468-009-9141-x}
\bibitem{BegMasLeray} E.J. Beggs and  I. Masmali, A Leray spectral sequence for noncommutative differential fibrations, Int. J. Geom. Methods Mod. Phys.  10 (2013) 1350015 (17pp.) \url{https://doi.org/10.1142/S0219887813500151}
\bibitem{Blake24} J.E. Blake, Submanifolds, Fibre Bundles, and Cofibrations in Noncommutative Differential Geometry, 2024 PhD thesis, \url{http://dx.doi.org/10.23889/SUThesis.67590}
\bibitem{BDMS} K. Bresser and F. M\"uller-Hoissen, A. Dimakis and  A. Sitarz, Noncommutative geometry of finite groups. J. Phys. A, 29 (1996) 2705--2735 \url{https://doi.org/10.1088/0305-4470/29/11/010}
\bibitem{conHig} A.\ Connes  and  N.\ Higson, D\'eformations, morphismes asymptotiques et K-th\'eorie bivariante, C.R.\ Acad.\ Sci.\ Paris S\'er.\  I Math.\ 311 (1990) 101--106
\bibitem{ConnesRieffel87} A. Connes and M. Rieffel, Yang–Mills for noncommutative two tori, Contemp. Math. 62 (1987) 237–266 \url{https://doi.org/10.1090/conm/062}
\bibitem{DVMass} M. Dubois-Violette and  T. Masson, On the first-order operators in bimodules, Lett.\ Math.\ Phys.\ 37 (1996) 467--474 \url{https://doi.org/10.1007/BF00312677}
\bibitem{DVMic} M. Dubois-Violette and  P.W. Michor, Connections on central bimodules in noncommutative differential geometry, J.\ Geom.\ Phys.\ 20 (1996) 218--232 \url{https://doi.org/10.1016/0393-0440(95)00057-7}
\bibitem{Echterhoff-fibrations} S. Echterhoff and R. Nest and H. Oyono-Oyono, Fibrations with noncommutative fibres, J. Noncommut. Geom. 3 (2009), no. 3, pp. 377–417 \url{https://doi.org/10.4171/jncg/41}
\bibitem{FioMad} G. Fiore and  J. Madore, Leibniz rules and reality conditions, Eur.\ Phys.\ J.\ C Part.\ Fields 17 (2000) 359--366 \url{https://doi.org/10.1007/s100520000470}
\bibitem{KaspKK} G. Kasparov. The operator K-functor and extensions of C*-algebras. Izv.\ Akad.\ Nauk.\ SSSR Ser.\ Mat.\ 44 (1980), 571-636
\bibitem{matrices-private-communication} T. Ivanova, personal communication, 2022
  \bibitem{Kaygun19} A. Kaygun, Noncommutative fibrations, Comm. Algebra 47 (2019), no. 8, 3384--3398, \url{https://doi.org/10.1080/00927872.2018.1559850}
\bibitem{Lance} E.C. Lance, Hilbert C*-modules, A toolkit for operator algebraists, LMS. Lecture Note Series 210, CUP. 1995 \url{https://doi.org/10.1112/S0024609396291672}
\bibitem{Madore} J. Madore, An introduction to noncommutative differential geometry and its physical applications, LMS Lecture Note Series, 257, CUP  1999 \url{https://doi.org/10.1017/CBO9780511569357}
\bibitem{Quantum-groups-calculus} S. Majid, Classification of bicovariant differential calculi, J. Geom. Phys. 25 (1998) 119–140 \url{https://doi.org/10.1016/S0393-0440(97)00025-9}
\bibitem{Majid16LTCC} S. Majid, Noncommutative differential geometry, in LTCC Lecture Notes Series: Analysis and Mathematical Physics, eds. S. Bullet, T. Fearn and F. Smith, World Sci. (2016) pp. 139– 176 \url{https://www.academia.edu/69863305/LTCC_Lectures_on_Noncommutative_Differential_Geometry}
\bibitem{S3-calculus} S. Majid and W-Q. Tao, Generalised noncommutative geometry on finite groups and Hopf quivers, J. Noncom. Geom. 13 (2019) 1055–1116 \url{https://doi.org/10.4171/jncg/345}
\bibitem{McCleary} J. McCleary, A User's Guide to Spectral Sequences, 2nd ed., Cambridge University Press, 2001 \url{https://doi.org/10.1017/CBO9780511626289}
\bibitem{Mourad} J. Mourad, Linear connections in noncommutative geometry, Class.\ Quant. Grav.\ 12 (1995)  965--974 \url{https://iopscience.iop.org/article/10.1088/0264-9381/12/4/007}
\bibitem{Spectral Sequence Image Wikipedia} L.penguu, Own work image, CC BY-SA 4.0, \url{https://commons.wikimedia.org/w/index.php?curid=110078656} last accessed Dec 2024
\bibitem{Spanier} E.H. Spanier, Algebraic Topology, Springer Science \& Business Media, 2012, ISBN 9781468493221 \url{https://doi.org/10.1007/978-1-4684-9322-1}
\end{thebibliography}
\end{document}